\documentclass[reqno,10pt]{amsart}

\usepackage[normalem]{ulem}
\usepackage{amsthm}
\usepackage[all]{xy}
\usepackage{graphics}
\usepackage[dvips]{graphicx}
\usepackage{amsmath}
\usepackage{amsfonts}
\usepackage{cancel}
\usepackage{rotating}
\usepackage{a4wide}
\usepackage{amssymb}
\usepackage{color}
\usepackage{fancybox}
\usepackage{wasysym}
\usepackage{textcomp}

\usepackage[english]{babel}
\usepackage[latin1]{inputenc}

\newcommand{\wt}{\widetilde}
\newcommand{\wh}{\widehat}
\newcommand{\Mg}{$(M,g)\,$}

\newcommand{\n}{\nabla}
\newcommand{\oln}{\overline{\nabla}}
\renewcommand{\O}{\Omega}
\newcommand{\ol}{\overline}
\newcommand{\ob}{\overbrace}
\newcommand{\beit}{\begin{itemize}}
\newcommand{\eit}{\end{itemize}}
\newcommand{\ben}{\begin{enumerate}}
\newcommand{\een}{\end{enumerate}}
\renewcommand{\t}{\theta}
\renewcommand{\o}{\omega}
\newcommand{\la}{\lambda}
\renewcommand{\a}{\alpha}
\renewcommand{\b}{\beta}
\newcommand{\s}{\sigma}
\renewcommand{\r}{\rho}
\newcommand{\g}{\gamma}
\newcommand{\G}{\Gamma}
\newcommand{\vp}{\varphi}
\newcommand{\p}{\partial}
\newcommand{\we}{\wedge}
\newcommand{\mc}{\mathcal}

\newcommand{\beq}{\begin{eqnarray}}
\newcommand{\eeq}{\end{eqnarray}}
\newcommand{\st}{Lorentzian manifold }

\newcommand{\tsy}{2nd-symmetric }

\newcommand{\tsns}{proper $2$nd-symmetric }
\newcommand{\nr}{partly null frame }

\def\R{\text{$\mathbb{R}$}}

\def\M{{\ol{\mc{M}}}}
\def\U{\mc{U}}

\def\K{{(T\mc{U})^\perp}}
\def\v#1{\vec{#1}}

\def\fr#1#2{\frac{#1}{#2}}

\def\mc#1{\mathcal{#1}}

\def\ni{\noindent}
\def\ri#1{\mathring{#1}}

\def\up{u_p}
\def\vp{v_p}
\def\xip{x^i_p}

\newtheorem{teorema}{Teorema}[section]

\newtheorem{teor}[teorema]{Theorem}
\newtheorem{claim}[teorema]{Claim}

\newtheorem{prop}[teorema]{Proposition}
\newtheorem{coro}[teorema]{Corollary}

\newtheorem{lema}[teorema]{Lemma}
\newtheorem{conve}[teorema]{Convention}
\newtheorem{rema}[teorema]{Remark}
\newtheorem{defn}[teorema]{Definition}

\title[Lorentzian manifolds with $\nabla^2R=0$]{Structure of second-order symmetric \\ Lorentzian manifolds} 
\author[OF Blanco]{Oihane F. Blanco}
\address{Departamento de Geometr\'{\i}a y Topolog\'{\i}a.
 Facultad de Ciencias, Universidad de Granada.
 Campus Fuentenueva s/n, E-18071 Granada, Spain}
\email{oihane@ugr.es}
\author[M. S\'anchez]{Miguel S\'anchez}
\address{Departamento de Geometr\'{\i}a y Topolog\'{\i}a.
 Facultad de Ciencias, Universidad de Granada.
 Campus Fuentenueva s/n, E-18071 Granada, Spain}
\email{sanchezm@ugr.es}
\thanks{The three authors are partially supported by the grant P09-FQM-4496 (J. Andaluc\'ia---FEDER), OFB and MS are partially
supported by MTM2010--18099 (Spanish MICINN---FEDER), and JMMS is
supported by grants FIS2010-15492 (MICINN) and GIU06/37 (UPV/EHU)}
\author[JMM Senovilla]{Jos\'e M. M. Senovilla}
\address{Departamento de F\'{\i}sica Te\'orica e Historia de la Ciencia. Facultad de Ciencia y Tecnolog\'ia, Universidad del Pa\'{\i}s Vasco, Apartado 644, E-48080 Bilbao, Spain}
\email{josemm.senovilla@ehu.es}


\begin{document}
\begin{abstract}
{\em Second-order  symmetric  Lorentzian spaces}, that is to say,
Lorentzian manifolds with vanishing second derivative $\nabla
\nabla R\equiv 0$ of the curvature tensor $R$,  are characterized
by several geometric properties, and  explicitly presented.
Locally, they are a product $M=M_1\times M_2$ where each factor is
uniquely determined as follows: $M_2$ is a Riemannian symmetric
space and $M_1$ is either a constant-curvature Lorentzian space or
a definite type of plane wave generalizing the Cahen-Wallach
family. In the proper case (i.e., $\nabla R \neq 0$ at some
point), the curvature tensor turns out to be described by some
local affine function which characterizes a globally defined
{parallel lightlike direction}. As a consequence, the
corresponding global classification is obtained, namely: any
complete second-order symmetric space admits as universal covering
such a product $M_1\times M_2$. From the technical point of view,
a direct analysis of  the second-symmetry partial differential
equations
 is carried out leading to several results of independent interest
 relative to spaces
 with a parallel lightlike vector field ---the so-called Brinkmann spaces.

\end{abstract}

\maketitle

\tableofcontents

\newpage

\section{Introduction}\label{s1}

A venerable result in Differential Geometry (Nomizu and Ozeki \cite{NO}, Tanno \cite{TA}
) states that, for a Riemannian manifold $(M,g)$, the vanishing
of the $r$-th covariant derivative of its curvature tensor $R$,
\begin{equation} \label{e00}
\n^r R\, (={\n\stackrel{(r)}{\dots}\n} R)\equiv 0 , \quad \quad r\geq 2,
\end{equation}
implies the vanishing of the first one,  i.e., that $(M,g)$ is
{\it locally symmetric}. As a consequence, the standard generalization of
Riemannian locally symmetric spaces are the 
semi-symmetric spaces, introduced by Cartan  \cite{CA2} and
defined by the commutativity of the covariant derivatives applied
to $R$: \beq\label{semi}R(X,Y)R:=
(\n_X\n_Y-\n_Y\n_X-\n_{[X,Y]})R=0 \quad \quad \text{ for all
vector fields } X,Y.\eeq Their structure was determined by Szab\'o
locally in \cite{SZ1} and globally later in \cite{SZ2}. However,
when $(M,g)$ is a Lorentzian manifold, the equality \eqref{e00}
does not imply $\nabla R=0$. In fact, the irreducible character of
de Rham decomposition, which was an essential ingredient in the
Riemannian result, fails for Lorentzian metrics. Thus, a ladder of
logical generalizations of Lorentzian locally symmetric spaces is
given by \eqref{e00}. We call these semi-Riemannian manifolds {\it
$r$th order symmetric} (or {\it $r$th-symmetric} for short)
spaces\footnote{$r$th-order symmetric spaces were introduced in
\cite{SN1} and termed as $r$-symmetric for short. However, a
different notion of {\em 3-symmetric space}, introduced by Gray
\cite{GR}, was already available and somehow spread in the
literature (for example, see the recent article \cite{GD}). Thus,
we have preferred to use the ordinal ({\em 3rd-symmetric}, say) to
avoid any possible confusion.}. The purpose of the present article
is to determine the simplest of these new classes explicitly: the
Lorentzian 2nd-symmetric spaces\footnote{The results of the
present paper were announced in particular in the Spanish
Relativity Meetings celebrated  in Bilbao, 7-11 September '09
(only the four dimensional case,
see \cite{BSS1}) and in
Granada, 6-10 September '10 where the general $n$-dimensional case was
considered (see \cite{BSS2}).}.

The classification of Riemannian locally 
symmetric spaces is known since Cartan's work \cite{CA} (see also
\cite{He,Be}), and the classification of the Lorentzian
simply-connected symmetric spaces was carried out by Cahen and
Wallach in \cite{CW}. Extensions to other signatures and to
non-simply-connected cases are also available, see Cahen and
Parker \cite{CP}, Neukirchner \cite{Neu} and specially Kath  and
Olbrich \cite{ka1,ka2}. Lorentzian semi-symmetric spaces have also
been studied in the literature, see for instance their
classification in four dimensions \cite{IE,SV} and references
therein. Nevertheless, prior to the paper \cite{SN1} by one of the
authors, the 2nd-symmetric spaces had not been studied
systematically. As pointed out in this reference, simple examples of {\em proper
$r$th-symmetric}
---$r$th-symmetric but not $(r-1)$th-symmetric--- Lorentzian
spaces can be constructed within the class of $n$-dimensional
plane waves, see subsection \ref{s3b} below. They constitute a
straightforward generalization of the locally symmetric {\em
Cahen-Wallach spaces} \cite{CW} and, as we will prove, they
essentially exhaust the whole class of \tsns Lorentzian spaces.

It is worth pointing out that the 2nd-symmetric spaces are
appealing from the viewpoint of the local group of symmetries of
the manifold, because the condition $\nabla^2 R=0$ can be
expressed in terms of the infinitesimal holonomy algebra of the
manifold (as this algebra is generated by the image of the
curvature two-form and its first derivative see, for example,
\cite[Th. 9.2 Ch. III]{KN}); in fact, the main result in
\cite{SN1} is a property of this holonomy group. With the help of
this property and the well-established results on locally
symmetric spaces by Cahen and Wallach, our proof will be
completely self-contained, by solving the equations of
2nd-symmetry crudely\footnote{While the present paper was being
finished, a different local approach fully based in holonomy
groups (including landmarks such as the classification of all the
holonomy groups in Lorentzian signature \cite{Le}) has been developed in \cite{AG},\cite{Gal}. By using it,
the crucial result  in \cite[Theorem 4.2]{SN1}
is revisited, and a partial version of the Theorem 1.1 below is provided.}.

Specifically, the main result we will prove is:
\begin{teor}\label{lem4}
 An  $n$-dimensional \tsns  Lorentzian space $(M,g)$ is
 locally isometric
 to a direct product $(M_1\times M_2,g_1\oplus
 g_2)$ where
 $(M_2,g_2)$ is a non-flat Riemannian symmetric space and $(M_1,g_1)$ is a {proper} generalized Cahen-Wallach
 space  of order $2$,
 defined as $M_1=\R^{d+2}$ ($d\geq0$) endowed with the metric
$$
g_1=-2du\left(dv+du\, \sum_{i,j=2}^{d+1} p_{ij}(u)x^ix^j
\right)+\sum_{i=2}^{d+1} (dx^i)^2,
$$
\noindent where $(u,v,x^2,\dots, x^{d+1})$ are the natural
coordinates of $\R^{d+2}$ and each function $p_{ij}$ is affine:
$p_{ij}(u)=\a_{ij} u+\b_{ij}$ for some $\a_{ij}, \b_{ij}\in \R$
{with al least one of the $\a_{ij}$ non vanishing,} and all
$i,j=2,\dots , d+1$.

Moreover, if $(M, g)$ is also geodesically complete and simply
connected, then $(M, g)$ is globally isometric to one such direct
product.
\end{teor}
Very roughly, the idea of the proof is the following. The starting
point is a   {significant result} obtained by one of the authors
(\cite[Theorem 4.2]{SN1}): {\it any simply-connected Lorentzian
\tsns space $(M,g)$ admits a parallel lightlike vector field $K$.}
Lorentzian spaces with such a $K$ were obtained by Brinkmann
\cite{BR} and will be studied in Sec.\ref{s100}, where local bases
associated to what we call {\em Brinkmann charts} $\{u,v,x^i\}$
will be introduced. The fact that, for any such chart, the slices
with constant $u$ and $v$ happen to be locally symmetric suggests
a reduction of the equations for $2$nd-symmetry. This reduction is
carried out by exploiting the integrability conditions in full,
and by applying some technical algebraic properties. Then,
Eisenhart-type decompositions can be used to transform the
original equations in $(M,g)$ into the 2nd-symmetry equations of
two
simpler spacetimes 
from which the stated parts $(M_1,g_1)$ and $(M_2,g_2)$ emerge
locally. Some geometric elaborations yield the global result. In
fact, along the proof it will become apparent that $(M,g)$ admits
a globally defined {parallel lightlike direction}. Thereupon, the
global requirements complete the proof easily.

Summing up, we will prove that the implicit symmetries in the
equations of 2nd-symmetry turn out to be sufficient to solve them,
and actually to find their general solution explicitly. Incidentally,
different technical tools for some classes of partial
differential equations, which may have interest in its own right,
will be developed.

\subsection{Outline of the paper}
This article is organized as
follows. In Section \ref{s2}, we fix all our conventions and explain the notation.
We are specially
careful with the latter, for we will use a powerful combination of
intrinsic expressions and tensor-component computations.

In Section \ref{s27}, we review some results on locally symmetric spaces (Subsection
\ref{s271}), and compare them with the known results on
2nd-symmetric ones (Subsection \ref{s272}). We also
describe the properties of the generalized Cahen-Wallach family of
Lorentzian manifolds (Subsection \ref{s3b}). They will turn out to be the
non-trivial part of the \tsns spaces (Proposition
 \ref{glob0}). Judicious interpretations alongside some technical results (Corollary \ref{lnuevo}, Lemma \ref{glob})
  will imply that the global counterpart of Theorem 1.1 can be obtained from the local one.
  Accordingly, in the following sections we will work locally, except when otherwise stated explicitly.

Section \ref{s100} is devoted to a local study of Brinkmann
spaces. It has a technical nature, but it may have interest on its
own right.
In Subsection \ref{s10} we revisit the
known procedure to find a {\em Brinkmann chart $\{u,v,x^i\}$} associated to a
{\em Brinkmman decomposition $\{u,v\}$}, so that the metric of the manifold is
written as:
\beq \label{m1}
g=-2du(dv+H(u,x^k) du+ 
W_i(u,x^k)
dx^i)+
g_{ij}(u,x^k) dx^idx^j,~~k=2,\ldots,n-1.
\eeq
(see Section \ref{s2} for notation). These yield an associated spacelike $(n-2)$-foliation $\M$ characterized by constant values of $u$ and $v$, and a non-necessarily
orthogonal timelike 2-foliation $\U$ generated by
$\partial_u, \partial_v$, as well as other distributions of
interest.
For the sake of clarity, the relations among tensor fields on
these distributions are briefly explained. In
Subsection \ref{difope}, three operators $\oln$, $\ol{d}$ and $\dot{}$
adapted to the foliations $\M$ and $\U$ are introduced and related to the connection $\n$, the exterior differential $d$ and the geometry of the foliation $\M$. In Subsection
\ref{s101}, computations {\em $\grave{a}$ la Cartan} on
local vector-field bases introduced in \eqref{frame}
yield manageable expressions of
the connection one-forms and the first two covariant derivatives of
the curvature tensor $R$. To help on simplifying the resulting expressions, a fourth differential operator $D_0$ with a transverse nature, which complements the three previous ones, is introduced. The dependence of all mentioned objects
on the Brinkmann chart and their behavior under changes of charts is emphasized and explicitly controlled. Then, a version of a classical Eisenhart theorem \cite{Eis}
adapted to our problem, which involves $u$-dependent metrics on the leaves of
a foliation endowed with the intrinsic $\oln$ and the
transverse $D_0$ derivatives, is provided in Theorem \ref{gei} of Subsection \ref{s7a}.
This allows us to obtain sufficient conditions such that the $u$-family of Riemannian metrics $\ol{g}=g_{ij}(u,x^k) dx^idx^j$ in \eqref{m1} are simultaneously reducible.

In Section \ref{s6} we solve the equations of $2$nd-symmetry on
Brinkmann spaces in several steps. Firstly, we prove that the equations of 2nd-symmetry imply that the foliation $\M$ must be locally symmetric for all Brinkmann charts $\{u,v,x^i\}$ associated to a fixed decomposition $\{u,v\}$, hence severe restrictions on
the curvature must hold. Secondly, using some auxiliary algebraic results on vector spaces, those restrictions are explicitly determined in Propositions
\ref{prop20} and \ref{prop21}. Moreover, a two-covariant tensor field $\tilde A$
on $\M$, associated to the Brinkmann decomposition $\{u,v\}$ only
but not to the other coordinates of the Brinkmann chart, is
defined in Corollary \ref{Aten}. All the information of the tensor field $\n R$ is
codified in $\tilde A$ and the leaves of $\M$. These results are summarized in Theorem
\ref{redu} and its corollary of Subsection \ref{s32}.
Thirdly, a reorganization of the $2$nd-symmetry equations into two independent blocks associated to different Brinkmann manifolds $(M^{[m]},g^{[m]})$ with $m\in\{1,2\}$
is proven in Subsection \ref{s8c} (Theorem \ref{cor2} and Proposition \ref{redBri}).
This is achieved by applying our version of Eisenhart theorem showing that the metric is
suitably reducible into a Ricci-flat part and a non-Ricci-flat one. Actually, the former is flat as a consequence of a known result \cite{Ale}
and the operator $\tilde A$ lives only
in this flat part (Proposition \ref{teo1}).
Finally, the proof of the main result is complete in Subsection \ref{s11}.
The first space $(M^{[1]},g^{[1]})$ is directly
computable leading to the required generalized Cahen-Wallach
expression of the part $(M_1,g_1)$ in Theorem \ref{lem4}. The
2nd-symmetry equations for the second space $(M^{[2]},g^{[2]})$
become equivalent to the equations for local symmetry so that $(M^{[2]},g^{[2]})$ collapses to a locally symmetric Brinkmann space, therefore the Cahen-Wallach classification allows us to determine the Riemannian locally symmetric part $(M_2,g_2)$ in Theorem
\ref{lem4}.

\bigskip

To end this introduction, we would like to emphasize that our
results open new questions and lines of interest. The first two
are obvious, consisting on the study of proper rth-symmetric
spaces with $r>2$, and the study of 2nd-order symmetric spaces
with metrics of index greater than 1. In both cases, our approach
is not directly applicable and, in fact, it is not clear how these
generalizations would affect even to our starting point
(\cite[Theorem 4.2]{SN1}). Moreover, recall that the obtained
proper 2nd-symmetric spaces share the symmetries of plane waves
(for explicit expressions including the more exotic Kerr-Schild
symmetries, see \cite{CHS}). So, an interesting question may be to
find connections between the group of symmetries inherent a priori
to rth-symmetry, and the actual symmetries of the eventually
obtained proper rth-symmetric spaces. Lastly, the non-simply
connected case and, in particular, the existence of compact
quotients of plane waves, becomes also a natural problem in this
setting (we acknowledge Professor A. Zeghib, from ENS Lyon, for
discussions stressing the importance of this question).

\section{Notation and conventions}\label{s2}

$M$ will denote a (connected) $n$-dimensional manifold. For simplicity, it will be implicitly assumed to be differentiable of class $C^k$ with
$k=\infty$, but one only needs $k=r+3$ for $r$-th
order symmetric spaces, i.e., $k=5$ for most of the paper. Accordingly, all objects will be assumed to be as differentiable as necessary depending on $k$.
For Lorentzian metrics, our convention on the signature is $(-,+,\ldots,+)$.
Indices written in greek small letters $\a,\b,\la,\ldots$
will run from $0$ to $n-1$, while those in latin small letters starting at $i$
$(i,j,k,\ldots)$ will run from $2$ to $n-1$ and the usual {\it summation convention}
is used.

A  chart on the
manifold will be indicated simply with its coordinate functions
$\{x^\a\}$.  When working on a Brinkmann space the coordinates
$x^0$ and $x^1$ will also be written as $u$ and $ v$, respectively, according to \eqref{m1}.
 For $f\in C^\infty(M)$ its partial derivatives are denoted by
\begin{equation}\label{e0}
\dot{f}\equiv\fr {\p f} {\p u},
\text{ and } f_{,i}\equiv\fr {\p f} {\p x^i}.
\end{equation}
Most of our computations will be local,
namely, in some appropriate neighborhood $U$ (reduced if neccesary) of any point
$p\in M$, but we will not specify the neighborhood ---as we have already done with the charts. In
general, we will use the notation as if $U=M$ except if there
were some possibility of confusion.

Let $TM$ and $T^*M$ be the tangent and the cotangent
bundles of $M$, respectively, and $\pi: T M\longrightarrow M$ the natural
projection.
A {\it $l$-distribution} of $M$ will be regarded as a $l$-subbundle $E$ of $TM$ and it will be involutive if $E=T\mc{F}$ for some foliation $\mc{F}$ of $M$. In this case, the bundle of all the $s$-covariant and
$r$-contravariant tensors on $T\mc{F}$ is denoted by $T^r_s\mc{F}$. The sections of a fiber bundle $\pi_E:
E\rightarrow M$ will be written as $\Gamma(E)$ unless the base is
not evident, in which case we use $\Gamma(M,E)$. In the case of sections of
$s$-forms the notation will be simplified:
$\Lambda^s E$ denotes the space of all the
$s$-form sections.

To write tensor equations in components, some local vector field
basis $\{V_{\a}\}$ on $TM$, or on some of its subbundles, plus its
dual basis $\{\zeta^\a\}$ for $T^*M$ are used. Of course, these
bases are not necessarily holonomic, i.e., associated to specific
coordinates $\{x^\a\}$, for which we use the standard notation
$\{\p_\a\},\{dx^\a\}$. We will follow typical notation for
covariant derivatives and their components as, for example, in
\cite[pp. 30-35]{Sa}. Notice that we denote the components of $R$
as defined in \eqref{semi} by $R^\a~_{\b\la\mu}$, which agrees
with \cite{HE,ES} but differs from \cite{Sa} where  the same is
written as  
${R_{\la\mu\b}~^\a}$ (therefore, $R_{\a\b\la\mu}$ differs on a
sign).

Sometimes the {\it abstract index notation} is also used, see \cite{Wal} for more information.
The {\it symmetrization}  (respectively {\it
anti-symmetrization}) of a tensor field $T$ is denoted by round (respectively square) brackets that enclose the indices to be symmetrized (respectively antisymmetrized). For example, $T_{(\a\b)\la} =(T_{\a\b\la}+T_{\b\a\la})/2$ and $T_{[\a|\b|\la]}=(T_{\a\b\la}-T_{\la\b\a})/2$. Furthermore, we write
$2dx^\alpha dx^\beta = dx^\alpha \otimes dx^\beta + dx^\beta
\otimes dx^\alpha$.
For the {\it wedge
product}, the convention  is
 $\b^{1}\wedge\ldots\wedge\b^{m}=\sum_{\sigma\in S_m}(-1)^{[\sigma]}\b^{\sigma(1)}\otimes\ldots\otimes\b^{\sigma(m)}$,
 where $\b^i$ are one-forms, $i=1,\ldots,m$ and $S_m$ denotes the set of all permutations of $\{1,\dots ,m\}$.
The metric isomorphism $\flat:T M \longrightarrow T^*M,
\v v\mapsto g(\v v,\cdot) $ and its inverse
$\sharp:T^*M\longrightarrow T M$, are written in components or abstract index notation so that
$X_{\a}:=(X^\flat)_\a= g_{\a\b}X^\b$ and
$\tau^{\a}:=(\tau^\sharp)^{\a}=g^{\a\b}\tau_\b$ for all $X\in\Gamma(TM),
\tau\in \Gamma(TM^*)$.

\section{Locally symmetric versus \tsy semi-Riemannian manifolds}\label{s27}
\subsection{Generalities on local symmetry}\label{s271}
{Let $(M,g)$ be a semi-Riemannian manifold  and $p\in M$. The {\em
local geodesic symmetry} $s_p$ with respect to $p$ is the
diffeomorphism $s_p:N_p\longrightarrow N_p$, defined on a
sufficiently small normal neighborhood $N_p$ of $p$, which maps
each $q=\g(1)\in N_p$ into  $s_p(q)=\g(-1)$, where $\g$ is the
univocally determined geodesic in $N_p$ from $p$ to $q$.} We
collect in the following proposition some characterizations of
locally symmetric semi-Riemannian manifolds, {i.e., manifolds
satisfying  $\n R=0$,} for comparison with \tsy spaces. See
\cite[pp.219-223]{ON2} for further details and results.
\begin{prop} \label{pls} For a semi-Riemannian manifold $(M,g)$ the following conditions are equivalent:
\ben \item [(i)] $(M,g)$ is locally symmetric.  \item[(ii)] If
$L:T_pM\longrightarrow T_qM$ is a linear isometry that preserves
curvature (i.e., for any $X_p,Y_p,Z_p\in T_pM$,
$L(R(X_p,Y_p)Z_p)=R(L(X_p),L(Y_p))L(Z_p)$) there exist small
normal neighborhoods $N_p$ of $p$ and $N_q$ of $q$ and a unique
isometry $\phi:N_p\rightarrow N_q$ such that $d\phi_{|p}=L$.
\item[(iii)] The local geodesic symmetry $s_p$ is an isometry at
any $p\in M$. \item[(iv)] If $X,Y$ and $Z$ are parallel vector
fields along a curve $\a$, then the vector field $R(X,Y)Z$ is also
parallel along $\a$. \item[(v)] The sectional curvature is
invariant under parallel translation: for any non-degenerate
plane, its parallel transport
along any curve has constant sectional curvature. \een
\end{prop}
Another important known result is \cite{He,ON2}:
\begin{prop}\label{pz}
 Any semi-Riemannian symmetric space is an analytic, (geodesically) complete,  homogeneous space $G/H$.
Moreover, the universal covering of a complete connected locally
symmetric space is symmetric.
\end{prop}
For the Riemannian case, we will need the following.
\begin{prop}\label{redlocsym}
 Let $(M,g)$ be a locally symmetric Riemannian manifold. Then,

\ben \item $(M,g)$ is locally isometric to the direct product of a finite number of
irreducible locally symmetric spaces and a Euclidean space of dimension $d\geq 0$.

\item If $(M,g)$ is irreducible then it is an Einstein manifold,
i.e., $\hbox{Ric}= c g$.

\item If $(M,g)$ is Ricci-flat (that is, Einstein with $c=0$), then it is flat.
\end{enumerate}
\end{prop}

{\em Proof.} (1) This a consequence of the classical de-Rham
decomposition of $M$, as any irreducible part must be locally
symmetric.

(2) As the Ricci tensor in a locally symmetric space is parallel,
the result follows from a classical result by Eisenhart (Theorem \ref{Eise1} below).

(3) By hypothesis, $(M,g)$ is locally isometric to a Ricci-flat
symmetric space and, by a result in \cite{Ale}, this space must be
flat\footnote{In fact, Alekseevskii and Kimelfeld \cite{Ale}
proved that any Ricci-flat homogeneous Riemannian space is flat.
It is worth pointing out that, as a difference with the locally
symmetric case,  the locally homogeneus spaces maybe non-regular,
that is, non-locally isometric to some homogeneous space
\cite{Ko}. However, Spiro \cite{Sp} showed that all locally
homogeneous spaces with non-positive Ricci curvature are regular,
hence the result in \cite{Ale} can be extended
to the locally homogeneous case. Nevertheless, it cannot be extended
to the Lorentzian case: it is easy to find a counterexample in the
Cahen-Wallach spaces below.}. $\square$

The classification of Lorentzian simply-connected symmetric spaces by Cahen and Wallach \cite{CW} can be
summarized as follows:
\begin{teor}\label{CW}
Any simply-connected Lorentzian symmetric space $(M,g)$ is isometric to the product of a simply-connected Riemannian symmetric space and one of the following
Lorentzian manifolds:
\begin{itemize}
 \item[(a)] $(\mathbb{R},-dt^2)$
\item[(b)] the universal cover of d-dimensional de Sitter or
anti-de Sitter spaces, $d\geq 2$, \item[(c)] A Cahen-Wallach space
$CW^d(A)=(\R^d,g_A)$, $d\geq 2$, where $A=(A_{ij})$ is a
$(d-2)\times (d-2)$ symmetric constant matrix, and the metric is
written \begin{equation}\label{edcw}g_A=-2du\left(dv+ A_{ij}x^ix^j
du\right)+ \delta_{ij} dx^i dx^j,\end{equation} where
$\delta_{ij}$ is the Kronecker delta (observe that $CW^2
=\mathbb{L}^2$, as $A$ necessarily vanishes).
\end{itemize}
Therefore, if a Lorentzian symmetric space admits a parallel lightlike vector field, then it is locally isometric to the product of a $d$-dimensional Cahen-Wallach
space and an $(n-d)$-dimensional Riemannian symmetric space with $d\geq 2$.
\end{teor}

\subsection{Characterization of \tsy spaces}\label{s272}
There is no a priori reason to think that the characterization for
2nd-symmetric Lorentzian manifolds may include properties on local
geodesic symmetries, as in Proposition \ref{pls}(iii), or any
other similar semi-local property. Of course, this does not mean
that such a property might {not} be found a posteriori.

For any geodesic $\gamma$ such that $p=\gamma(0)$ and tangent plane $\Pi\subset T_pM$, consider the parallelly transported plane $\tau
\mapsto \Pi_\gamma(\tau) \subset T_{\gamma(\tau)}M$. Using
that the scalar product between parallelly propagated vector fields is constant
and that the sectional curvature $K(\Pi_p)$ on non-degenerate planes $\Pi_p$ determines
the full curvature tensor, one easily derives the next Lemma.

\begin{lema}\label{seccur}
The following three conditions are equivalent:
\ben \item[(i)] For any non-degenerate tangent plane
$\Pi$, its parallel transport $\Pi_\gamma$ {along any geodesic
$\gamma$} satisfies that $\fr d {d\tau} (K(\Pi_\g))$ remains
constant along $\g$. \item[(ii)] For any parallelly propagated
vector fields $X, Y, Z$ along {any geodesic $\gamma$}, the vector
field $(\n_{\g'} R)(X,Y)Z$ is itself parallelly propagated along
$\gamma$.
\item[(iii)] $\n_X(\n_Y R)-\n_{\n_X Y}R$ is skew-symmetric in $X,Y$.
\een

Moreover, if these  conditions hold, then the following property
follows:

\ben \item [(S)] if $\Pi$ is a lightlike plane with radical
spanned by $\v v\in T_pM$, and $V$ and $\Pi_{\g}$ are the parallel
transports of $\v v $ and $\Pi$ along {any geodesic} $\g$,
respectively, then $\fr d {d\tau} (K_{V}(\Pi_{\g}))$ remains
constant along $\g$, where  $K_{V}(\Pi_{\g})$ denotes the {\em
null sectional curvature along $\g$} and is defined by \cite{Ha}
\beq\nonumber K_{V}(\Pi_\g)=\fr {R(V,X,V,X)} {g_p(X,X)} \eeq with
$\{V,X\}$ spanning $\Pi_\g$ and $X$ parallel along $\g$.\een
\end{lema}
This result leads to the sought characterization of 2nd-symmetry.
\begin{prop}\label{2sdef}
The following statements are equivalent for a Lorentzian manifold \Mg:
 \ben
\item[(i)] $\n\n R=0$, i.e., \Mg \, is \tsy. 
\item[(ii)] If $V,X,Y,Z$ are parallelly propagated vector fields
along any curve $\a$, then $(\n_V R)(X,Y)Z$ is itself parallelly
propagated along the curve. \item[(iii)] $(M,g)$ is semi-symmetric
{(i.e., the curvature tensor fulfills \eqref{semi})} and satisfies
the  equivalent conditions in Lemma \ref{seccur} for any geodesic
$\g$. \een
\end{prop}

{\em Proof}. As (ii) characterizes when the tensor $\n R$ is
parallel, this condition is equivalent to (i). Also, both
conditions imply (iii) trivially. For the converse, notice that,
by semi-symmetry, ${\n\n R(\cdot~;U_1,U_2)=}\n_{U_2}(\n_{U_1}
R)-\n_{\n_{U_2} U_1}R$ is symmetric in $U_1, U_2$, 
and, {\it a fortiori}, it vanishes by applying the condition (iii) of
Lemma \ref{seccur}. 
 $\square$

\begin{rema}
 {\rm For any  Lorentzian manifold that satisfies the
equivalent conditions in Lemma \ref{seccur}
the (null) sectional curvature
of a plane parallelly propagated
along a geodesic varies as an affine function on the affine
parameter of the geodesics; in other words, the curvature along geodesics grows linearly.
Recall that the constancy of the sectional curvature characterized
locally symmetric spaces (Proposition \ref{pls}). On the other hand, there are obvious situations in which its non-constant linear growth
can be excluded. For
example, assume that $\gamma:[0,\infty) \rightarrow M$ is a
complete lightlike half-geodesic such that its velocity is
imprisoned in a compact subset $C\subset TM$. Then, for any
lightlike plane $\Pi$ with radical spanned by $\v v=\gamma'(0)$,
the derivative of $K_{V}(\Pi_{\g})$ must vanish. (There exists a sequence $\tau_n \nearrow
\infty$ such that $\{\gamma(\tau _n)\} \rightarrow q\in M$,
$\{\gamma '(\tau _n)\}\rightarrow \v x_q \in T_{q}M$ and
$\{\Pi_\gamma (\tau_n)\}\rightarrow \Pi_q \subset T_{q}M$, hence
$\{K_{V}(\Pi_{\g}(\tau_n))\}\rightarrow K_{\v x_q}(\Pi_q)$ and
this is incompatible with a linear growth of the curvature.)
This sort of properties, together with Proposition \ref{2sdef}, suggests global obstructions to the
existence of compact \tsns spaces. }
\end{rema}

\begin{coro} \label{lnuevo}  {Let $(M,g)$ be  a proper
(connected) 2nd-symmetric space. Then, $\n R\neq 0$ holds
everywhere, and $(M,g)$ admits a unique {parallel lightlike
direction}.}

\end{coro}

\noindent {\em Proof.} Obviously, if  $(\nabla R)_p\neq 0$ at some
$p\in M$, the (parallel tensor) $\nabla R$ cannot vanish at any
point. To prove the existence of the {parallel lightlike
direction} recall that, as we have already mentioned, around each
point $p$ there exists a parallel lightlike vector $K$
\cite[Theorem 4.2]{SN1}. Moreover, there cannot exist a second
such $K'$ that is independent of $K$ at $p$, as otherwise a
parallel timelike vector field $T$ could be constructed as a
linear combination of $K$ and $K'$; hence, the metric would split
around $p$ as a product $-dt^2\oplus g_R$, with $T=\nabla t$ and
$g_R$ a Riemannian metric, which should be locally symmetric, in
contradiction with $(\nabla R)_p\neq0$. Thus, the corresponding
parallel lightlike directions locally generated by $K$ and $K'$
must agree, so that they match in a single global one. $\Box$

\subsection{Generalized Cahen-Wallach spaces of order $r$}\label{s3b}
In this section we introduce the archetypes for $r$th-symmetric Lorentzian
manifolds, which are generalizations of the $d$-dimensional Cahen-Wallach spaces $CW^d(A)$
 introduced  in Theorem \ref{CW}.

Consider the larger family $CW^d_r(A)$ of all the
 {\em $d$-dimensional generalized Cahen-Wallach spaces of order $r$} $CW^d_r(A)=(\R^d,g_A)$, defined by using the same expression \eqref{edcw} but letting now $A$ to
depend on $u$ as a matrix of polynomials of degree less or
equal than $r-1$
\beq\label{ecw}A\equiv A(u)=A^{(r-1)}u^{r-1}+\ldots+A^{(1)} u+ A^{(0)}\eeq
\ni where each $A^{(l)}$ is a constant symmetric $(d-2)\times (d-2)$ matrix
for all $l\in\{0,\dots , r-1\}$. As before, a
generalized Cahen-Wallach space is {\it proper} if at least one of the polynomials has degree $r-1$, that is, if $A^{(r-1)}\neq 0$.
Notice that all the spaces $CW^d_r(A)$ contain a parallel
lightlike vector field (see Section \ref{s10}), and that $CW^d_1(A)=CW^d(A)$.
 In addition
\begin{prop}\label{glob0} Any proper generalized Cahen-Wallach space
is analytic, (geodesically) complete and proper $r$th-symmetric.
\end{prop}

{\it Proof}. The analyticity is obvious and the completeness follows from \cite[Proposition 3.5]{CFS}, where a more general type of plane waves was treated.
The last property was already mentioned in \cite{SN1} and follows by
computing the derivatives of the curvature
tensor  in the basis $\{E_\a\}=\{\p_u-H\p_v,\p_v,\p_i-W_i \p_v\}$
(see Section \ref{s101}). 
The only non-vanishing
components of $\n^l R$
for  $l\in\{0,\dots , r-1\}$ (for $l\geq r$, $\n^l R=0$) are
$$\stackrel{(l)}{\n_0\ldots\n_0} R^{1}~_{i0j}=\fr {d^l A_{ij}}{d u^l} =\sum_{k=l}^{r-1} \frac{k!}{(k-l)!} A_{ij}^{(k)} u^{k-l}$$
which leads to the result immediately.
$\square$

The following  lemma will be used to reduce the global version of
Theorem \ref{lem4} to the local one.
\begin{lema}\label{glob}
Let $(M,g)$ be a complete simply-connected Lorentzian manifold
which is locally isometric to the product of some
generalized Cahen-Wallach space with a Riemannian symmetric
space. Then, $(M,g)$ is in fact globally isometric to one of such
products.
\end{lema}

{\em Proof.} By assumption, $(M,g)$ is locally isometric to an
analytical manifold due to Propositions \ref{pz} and \ref{glob0}
and, thus, it is analytical too. 
The result follows from the
fact that, for any two complete simply-connected analytic
semi-Riemannian manifolds $(M,g)$ and $(M',g')$, every isometry
defined between connected open subsets of $M$ and $M'$ can be
uniquely extended to an isometry for the entire $M$ and $M'$ (see, for
example \cite[Cor. 6.4 Ch. VI]{KN} for the Riemannian case, and \cite[Th.
6.1 Ch. VI]{KN}, \cite[Cor. 7.29]{ON2} for
its generalization to the semi-Riemannian one). $\square$

\section{Brinkmann spaces}\label{s100}
A Lorentzian manifold is called a {\em Brinkmann space} if and only if it admits a parallel lightlike vector field. Brinkmann spaces have deserved increasing attention in recent years, see e.g. \cite{Ba},\cite{BL},\cite{FI}.
In this section, we derive some of their properties relevant to our problem.
\subsection{Basics on Brinkmann spaces}\label{s10}
In what follows, the lightlike parallel vector field $K$ of a Brinkmann space $(M,g)$ will always be assumed to be fixed. Then, it is well-known that
 any point $p\in M$ admits a coordinate chart
$\{x^\a\}=\{u,v,x^i\}$, that we call a {\it Brinkmann chart}, in such a way
that the metric takes the form \eqref{m1} 
and $K=-\partial_v$ \cite{BR}, see also \cite{ZA}. Without loss of
generality, we will assume that the range of the coordinates
includes $|u|, |v|, |x^i| <\epsilon$ for some $\epsilon>0$. 

\begin{figure*}
  \centering
    \includegraphics[width=0.7\textwidth]{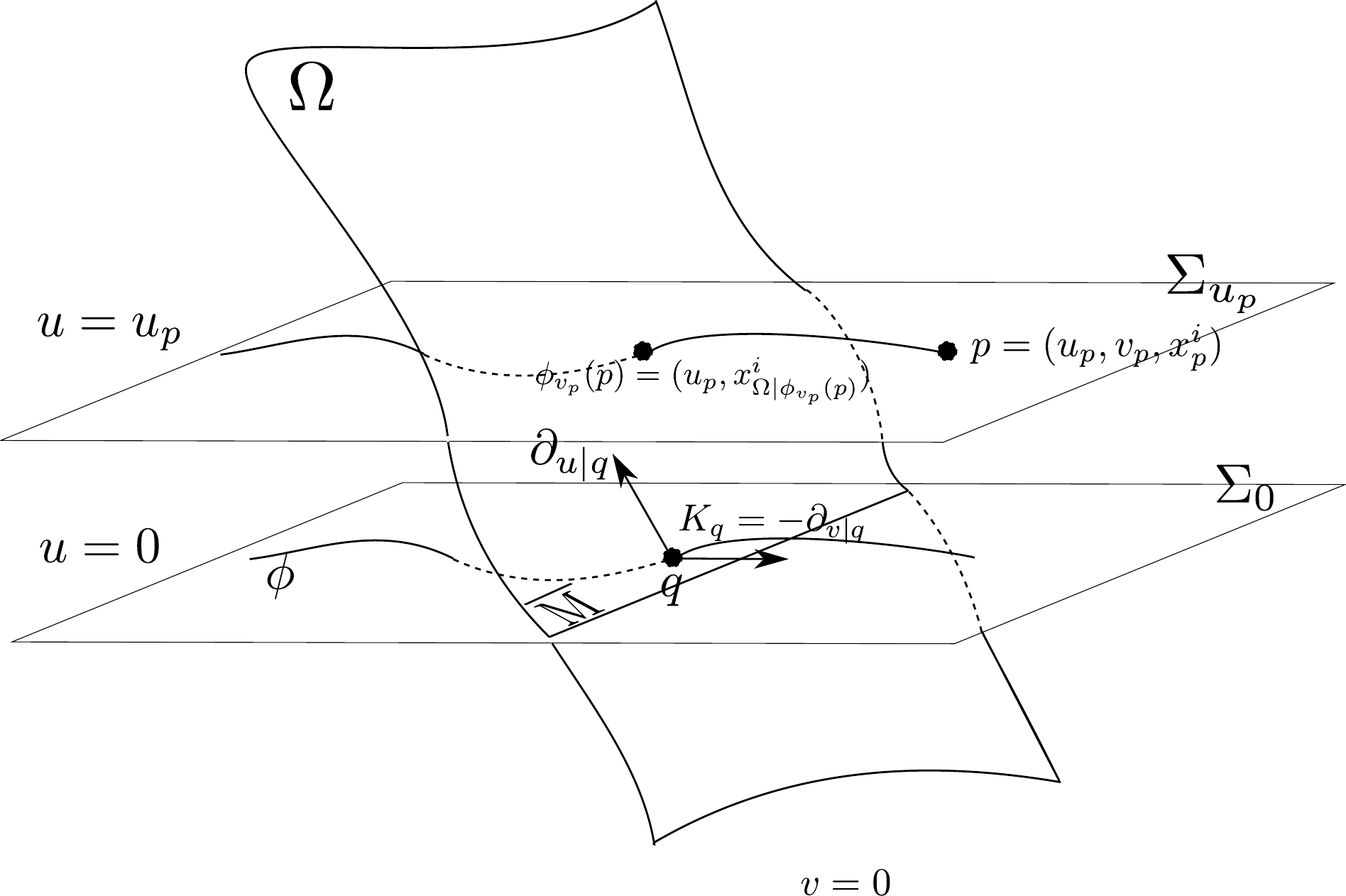}
\caption{\footnotesize Construction of a Brinkmann chart for a fixed lightlike parallel vector field $K=-\p_v$.}                                         
  \label{fig1}
\end{figure*}

Brinkmann charts can be obtained by means of the following
process (see Figure \ref{fig1}).
 As $K$ is parallel, choose a function $u$ such that
$K=\nabla u$ and with the value $0$ in the image of $u$. Each level
set $\Sigma_{u_0}=u^{-1}(u_0)$ is a lightlike integral
manifold of the distribution $K^\perp$ orthogonal to $K$. Set 
$\Sigma=\Sigma_0$ and choose a hypersurface $\Omega$
which is transverse to both $\Sigma$ and $K$. The function $u$ will serve as a coordinate
for $M$ as well as for $\Omega$. In $\Omega$ we choose a coordinate neighborhood
$\{u,x^i_{\Omega}\}$ completing $u$. The coordinate $v$ on $M$ is defined
by using the flow $\phi$ of $K$ to move each point $p$
to  
the point $\phi_{v_p}(p)\in \Sigma_{u_p}\cap\Omega $, and then we put $x^i_p =
x^i_{\Omega}(\phi_{v_p}(p))$. Observe that we have chosen to move
$p$ by using  the flow of $K=-\p_v$ so that $1\equiv
du(\p_u)=g(\nabla u,\p_u)=g(K,\p_u)=-g(\p_v,\p_u)$.

Conversely, the expression \eqref{m1} selects $K$ as
$-\partial_v$, and $\Sigma$ and $\Omega$ are the hypersurfaces $u=0$ and
$v=0$ respectively. Locally, a pair $(\Sigma, \Omega)$ determines a
{\em Brinkmann decomposition}, i.e., a pair of functions
$\{u,v\}$  constructed as above, which may serve as the two first
coordinates for different Brinkmann charts.

If  $\{u',v',x'^i\}$ denotes a second Brinkmann chart which
overlaps $\{u,v,x^i\}$, the corresponding $\Sigma'$ will also be a
level set of $u$, and $\Omega'$ can be regarded as a graph on
$\Omega$. So, the  change of coordinates can be written
as
 \beq
\label{chcoo}u'=u-u_0,\hspace{0.5cm} v'=v+F(u,x^j);\hspace{0.5cm}
x'^i=x'^i(u,x^j) \eeq 
Consequently, the relations between $H$, $W_i$, ${g}_{ij}$ in the original coordinates $\{u,v,x^i\}$ and $H'$, $W'_i$, ${g}'_{ij}$ in the new ones $\{u',v',x'^i\}$ are:
\beq
\label{H}
H&=&H'+\dot{F}+W'_{i}\dot{x}'^{i}-\frac{1}{2}{g}'_{ij}\dot{x}'^{i}\dot{x}'^{j}\\
\label{Wi}
W_i&=&W'_{j}\frac{\partial x'^{j}}{\partial x^{i}}+F_{,i}-\frac{1}{2}{g}'_{jk}\left(\frac{\partial x'^{j}}{\partial x^i}\dot{x}'^{k}+\frac{\partial x'^{k}}{\partial x^i}\dot{x}'^{j}\right)\\
\nonumber \label{gG}
g_{ij}&=&{g}'_{kl}\frac{\partial {x'}^{k}}{\partial
x^i}\frac{\partial {x'}^{l}}{\partial x^j} \eeq

The freedom in the choice of $\Omega$ makes it possible to obtain a
Brinkmann decomposition $\{u,v\}$ with a Brinkmann chart such that $H\equiv 0 \equiv W_i$, in
particular $\partial_u$ is lightlike and geodesic in the associated chart,
see e.g.
\cite{S}. For the sake of completeness, let us construct such a coordinate chart. Choose some $\Sigma$ as above and take any $(n-2)$-submanifold
$\ol{M} \hookrightarrow \Sigma$ which is transverse to $K$.
All such $\ol{M}$ are locally isometric, as $K$ is a
(parallel) lightlike direction in $\Sigma$. Now, consider for each $x\in
\ol{M}$ the {\em unique lightlike direction} $\v l_x$ orthogonal
to $\ol{M}$ and linearly independent of $K_x$. In a
small neighborhood, construct $\Omega$ by taking the geodesics with initial
velocity $\v l_x$ for all $x\in \ol{M}$. Complete the chart by choosing some local
coordinates $\{x^i_{\ol{M}}\}$ in $\ol{M}$ and defining the
coordinates $x^i_\Omega$ at each $y\in \Omega$ by $x^i_\Omega
(y):= x^i_{\ol{M}}(x_y)$, where $x_y$ is the unique point in
$\ol{M}$ which lies on the same lightlike geodesic as $y$.
Notice that $\Omega$ is a lightlike hypersurface, and the
corresponding coordinate vector field $\partial_u$ spans  its
radical, i.e., $H=0=W_i$, as required. We emphasize that the value of
$H$ and $W_idx^i$ depend on the choice of the coordinates $\{x^i\}$.
Moreover, given a Brinkmann decomposition $\{u,v\}$, there exists
a Brinkmann chart $\{u,v,x^i\}$ with $H=0=W_i$ if and only if the hypersurface
$\Omega$ obtained as $v=0$ is lightlike (when
$\Omega$ is lightlike, the integral curves of its radical must be
lightlike pregeodesics of $M$, because of the local maximizing
properties of these curves).
In spite of its simplicity, such a Brinkmann decomposition will
not be especially relevant for our study. It might simplify some
intermediate computations, but they are not well adapted to the generalized
Cahen-Wallach spaces of order $2$ 
which will turn out to be, as already announced, the essential part of proper 2nd-symmetric Lorentzian manifolds.

Given a Brinkmann chart, each integral curve of  $\p_v$ is
labelled by $\{u=u_0,x^i=x^i_0\}$, and each hypersurface
$\Omega_{v_0} \equiv \{v=v_0\}$ is a general
pseudo-Riemannian hypersurface (possibly signature-changing, as studied
systematically in \cite{SM}) and isometric to $\Omega$. 
We will repeatedly use the two natural transverse foliations
associated to each Brinkmann decomposition, namely:

\ben \item The $(n-2)$-dimensional foliation $\M$ with leaves
$\ol{M}=\Sigma\cap \Omega$ and the submanifolds obtained by moving $\ol{M}$ with the
flows of $\partial_u$ and $\partial_v$. Each leaf is defined by
$\{u=u_0,v=v_0\}$ and represented by
$\overline{M}_{(u_0,v_0)}$. The induced metric will be denoted by
$\ol{g}$ ($\ol{g}: T\M \times T\M \rightarrow \R$), so that
$\ol{g}_{ij}=g_{ij}$ in any Brinkmann chart. When necessary, $\ol{g}$ will be regarded as
a metric on a single leaf.
This foliation depends only on the Brinkmann decomposition $(\Sigma, \Omega)$, or equivalently on the chosen functions $\{u,v\}$ only.

\item The 2-dimensional foliation $\U$ whose leaves are the surfaces
obtained by moving each single point with the flows of
$\partial_u$ and $\partial_v$. Each leaf is given by $\{x^i=c^i_0\}$ for some constants $c^i_0$
and the induced metric is 
$-2du(dv+Hdu)$. This foliation depends on the expression of
$\partial_u$ and, thereby, on the coordinates chosen on $\Omega$ for
the Brinkmann chart. \een

Any Brinkmann chart  allows for the decomposition of the tensor bundles
in different ways, some of them to be detailed here. By a {\it \nr}  $\{E_\a\}$ for a \st we  mean a local
basis of  vector fields  such that
\begin{center}
\begin{tabular}{lll}
 $g(E_0,E_0)=0;$&$g(E_1,E_1)=0;$&$g(E_0,E_1)=-1;$\\
$g(E_0,E_i)=0;$&$g(E_1,E_i)=0;$&$g(E_i,E_j)=a_{ij},$
\end{tabular}
\end{center}
\ni for some functions $a_{ij}=a_{ji}$, which must obviously define a positive
definite metric. Its dual basis will be denoted by $\{\theta^\a\}$.
Now, fix a Brinkmann chart and  consider the associated
foliations $\M$, $\U$ and distributions $T\M$, $T\U$ and $\K$. The (canonical) 
{\em {\nr}  of the Brinkmann chart} is given by 
\begin{equation}
\{E_\a\}=\{\p_u-H \p_v,\, \p_v,\, -W_i\p_v+\p_i\},\,
\{\theta^\a\}=\{du, dv+H du+W_j dx^j, dx^i\},~~~ i=2,\ldots,n-1
\label{frame} \end{equation} 
which has $a_{ij}=g_{ij}$. Notice that $\K$ is spanned by the vector
fields $\{E_i\}$, and the brackets $[E_i,E_j]=(\p_j W_i-\p_i
W_j)\p_v \in \Gamma(T\U )$ measure its lack of involutivity.

The decomposition $TM= T\U \oplus T\M$ associated to any Brinkmann chart yields, on the one hand, the projection
$$\mathcal{P}_{\M}: \G(TM)\longrightarrow \G(T\M), \, \, X\mapsto\ol{X}$$
so that $\mathcal{P}_{\M}(\p_u)=\mathcal{P}_{\M}(\p_v)=0$ and
 $\text{$\mathcal{P}$}_{\M}(E_i)=\p_i$, and on the other hand the natural inclusion $$\mathcal{I}_{\M}:\G(T\M)\longrightarrow \G(TM)$$ with the induced dual map by restriction 
$$\text{$\mathcal{I}$}^*_{\M}:\G(T^* M)\longrightarrow\G(T^* \M), \,\, \b\mapsto\ol{\b}.$$
Clearly,
$\{\ol{E_i}\}=\{\p_i\}$ is a basis in $T\M$ whose cobasis is $\{\ol{\t}^i=\ol{dx^i}\}$. If
$X=X^\a E_\a\in \G(TM)$ and $\b=\b_\a \t^\a\in \G(T^* M)$,
then $\ol{X}^i=X^i$ and
$\ol{\b}_i=\b_i$ in the given bases.

We also introduce two linear homomorphisms between spaces of sections defined as:

(1) $\bar{ } :\G(T^r_s M)\longrightarrow\G(T^r_s\M)$, which maps each $T=T^{\a_1\ldots \a_r}_{\b_1 \ldots \b_s}{E}_{\a_1}\otimes\ldots\otimes{E}_{\a_r}\otimes{\t}^{\b_1}\otimes\ldots\otimes{\t}^{\b_s}$ into $ \ol{T}=T^{i_1\ldots i_r}_{j_1 \ldots j_s}\ol{E}_{i_1}\otimes\ldots\otimes\ol{E}_{i_r}\otimes\ol{\t}^{j_1}\otimes\ldots\otimes\ol{\t}^{j_s}$,

(2) $\ri{ } : \G(T^r_s \M)\longrightarrow\G(T^r_s M)$, which maps each $T=T^{i_1\ldots i_r}_{j_1 \ldots
j_s}\ol{E}_{i_1}\otimes\ldots\otimes\ol{E}_{i_r}\otimes\ol{\t}^{j_1}\otimes\ldots\otimes\ol{\t}^{j_s}$ into $\ri T=T^{i_1\ldots i_r}_{j_1 \ldots
j_s}{E}_{i_1}\otimes\ldots\otimes{E}_{i_r}\otimes{\t}^{j_1}\otimes\ldots\otimes{\t}^{j_s}$.

In general, $\ol{\ri T}=T$ but $\ri{\ol{T}}\neq T$,
because in the first case $T\in\G(T^r_s \M)$ necessarily,
but not in the second.

\begin{rema}\label{propied}{\rm  The previous maps are $C^\infty(M)$ linear
homomorphisms between spaces of sections, they commute with
contractions,
tensor products and, when applicable, wedge products in a natural
way; moreover, they leave invariant the $C^\infty(M)$ functions
considered as (0,0)-tensor fields: $\ri f=\ol{f}=f$ for all $f\in C^\infty(M)$.
Both homomorphisms $\bar{ }$ and $\ri{ }$ have trivial expressions in the introduced
bases. Essentially, all components of the
tensor fields and their images remain equal except for the fact
that, in the case of $T\mapsto\ol{T}$, the components aligned with any of $E_0$, $E_1$, $\theta^0$ or $\theta^1$ must be dropped while, in the case of $T\mapsto\ri T$, these components must
be restored with vanishing value. This simplifies our subsequent work in components substantially, as we will not need to distinguish notationally among different tensor fields derived
from a single one: it will be enough to realize which space of sections
is being considered. }
\end{rema}
As important illustrative examples, observe that the functions $W_i$ in the Brikmann expression \eqref{m1} can be regarded as the components in the basis
$\{\ol{E_i}\}=\{\partial_i\}$ of a one-form section in $\M$, namely: $W=W_i \ol{\t}^i\in\G(T^* \M)$, as well as the components in the basis
$\{{E_\a}\}$ of the one-form section $\ri W$ in $M$. Analogously, ${g}_{ij}$
can be regarded as the components 
of the projection  $P_\M(g)=\ol{g}\in\G(T^0_2 \M)$, which is the
inherited metric on $\M$. Therefore, the metric $g$ can be rewritten as
$$g=-2du(dv+Hdu+\ri W)+\ri {\ol{g}}.$$
One should also keep in mind that, when a single leaf $\ol{M}_{(u,v)}$ is
considered, ${g}_{ij}$ will also denote the components of the
induced metric $\ol{g}$ on the leaf, with no explicit mention to $(u,v)$ or the underlying Brinkmann decomposition.

\subsection{Three differential operators adapted to $\M$ and $v$-invariant sections}\label{difope}
By using the previously obtained vector bundle decomposition
associated to each Brinkmann chart, three differential operators
$\ol{d}$, $\oln$ and $\dot{} $ (``dot'') are introduced next. The operator $\oln$ depends on the Brinkmann decomposition
$\{u,v\}$ only, while the other two depend on the whole Brinkmann chart $\{u,v,x^i\}$.

\subsubsection{$v$-invariant tensor fields}

The variation on $\M$ of $T\in\G(T^r_s \M)$ under displacements
along the direction $E_1=\partial_v$
will be rather irrelevant for us. In fact, $E_1=-K$ is parallel and,
therefore, all the interesting objects to be  used here
 will also be invariant by its flow. More precisely:
\begin{defn}\label{isom}
 A tensor field $T\in\G(T^r_s \M)$ on $\M$ is {\em $v$-invariant} if $\ri T$ is  Lie-parallel along the flow of $E_1$, that is,
 $\mc{L}_{E_1} \ri T=0$.
\end{defn}

\begin{rema}\label{rlata} {\em
(1)  As $E_1=\p_v$, a tensor field $T=T^{i_1\ldots i_r}_{j_1\ldots
j_s}\p_{i_1}\otimes\ldots\otimes\p_{i_r}\otimes\ol{d}x^{j_1}\otimes\ldots\otimes\ol{d}x^{j_s}\in\G(T^r_s
\M)$ is $v$-invariant  if and only if $\p_v(T^{i_1\ldots i_r}_{j_1\ldots
j_s})=0$. For example,  $\ol{g}$, $H$ and $W$
are $v$-invariant.

(2) Obviously, if $T$ is $v$-invariant then it is determined by
its value on any of the hypersurfaces $\Omega_{v_0}=\{v=v_0\}$, and any
section $T_{\Omega_{v_0}} \in\Gamma(\Omega_{v_0}, T^r_s\M)$, defined only on
$\Omega_{v_0}$, can be extended to a unique $v$-invariant section
$T_\M\in \Gamma(T^r_s\M)$.

(3) As it will become obvious later, the derivatives of a
$v$-invariant section with respect to any of the three operators $\ol{d}$, $\oln$ and $\dot{ }$ to be defined below  are also $v$-invariant. Hence, these operators
can be naturally defined for sections in
$\Gamma(\Omega_{v_0}, T^r_s\M)$ (for instance, $\oln T_{\Omega_{v_0}}$ is the restriction
of $\oln T_\M$ to $\Omega_{v_0}$).
 }\end{rema}

\ni There is an alternative definition of $v$-invariance.
\begin{prop}\label{corE1}
$T\in\G(T^r_s \M)$ is $v$-invariant if and only if $\ri T$ is parallel in the
direction $E_1$.
\end{prop}
\ni {\it Proof.} As $E_1$ is  parallel, $\n_{E_1}
Q=\mc{L}_{E_1} Q$ for any $Q\in \G({T}^r_s M)$. $\square$

\subsubsection{The $\M$ exterior derivative $\ol{d}$}
\begin{defn} The {\em
$\M$ exterior derivative} $\ol{d}: \Lambda^{s}\M \rightarrow
\Lambda^{s+1}\M$ associated to  a Brinkmann chart
 $\{u,v,x^i\}$ is defined by:
 $$\ol{d}\b=\ol{d\ri \b} 
 \quad \quad
 \forall\b\in\Lambda^s \M,$$
where $d$ is the usual exterior derivative on the manifold $M$.
\end{defn}
It is straighforward to check that, if $\b\in\Lambda^s M$, then $\ol{d}~\ol{\b}=\ol{d\b}$.
In particular, $\ol{d}f=\ol{df}$ for $f\in C^\infty(M)$ and
thus $\ol{\t^i}=\ol{dx^i}=\ol{d}x^i$. This allows us to use the expressions $\p_i$ and
$\ol{d}x^i$ instead of $\ol{E_i}$ and $\ol{\t}^i$,
respectively.  Moreover, using the notation introduced in
\eqref{e0}, \beq\nonumber
df=\dot{f}du+\p_v f dv +\ri
{\overbrace{\ol{d}f}}.\eeq  

\begin{prop}\label{df1}
 Let $\b\in\Lambda^s\M$. The differential $\ol{d}$ satisfies the following properties:
\ben \item It is linear and
$\ol{d}(\b\wedge\tau)=\ol{d}\b\wedge\tau+(-1)^s
\b\wedge\ol{d}\tau$ for all $\tau\in \Lambda^q\M$.
 \item $\ol{d}(\ol{d}\b)= 0$.
\item 
If $\b=\fr 1 {s!}\b_{i_1\ldots
i_s}\ol{d}x^{i_1}\wedge\ldots\ol{d}x^{i_s}$, then $\ol{d}\b=\fr 1
{s!}\p_k(\b_{i_1\ldots i_s})
\ol{d}x^k\wedge\ol{d}x^{i_1}\wedge\ldots\ol{d}x^{i_s}$
 \item
(Poincar\'e Lemma).
If $\ol{d}\beta =0$, then each $p\in M$ admits a neighborhood $U$
and a section of $(s-1)$-forms $\tau$ such that $\ol{d}\tau
=\beta$ on $U$. Moreover, if the
$\ol{d}$-closed $s$-form $\beta$ is $v$-invariant, then  so can be
chosen the  $(s-1)$-form $\tau$.
\een
\end{prop}

\ni {\em Proof.} (1) Straightforward.

(2) $\ol{d}(\ol{d}\b):=\ol{d}(\ol{d\ri \b})=\ol{d(d\ri \b) }=0$.

(3) Apply (1), (2) and $\ol{d}f=\ol{df}$ in $\ol{d}(\fr 1
{s!}\b_{i_1\ldots i_s}\ol{d}x^{i_1}\wedge\ldots\ol{d}x^{i_s})$.

(4) We sketch the  notationally simpler case $s=1$ (to be used
later). For $\beta=\beta_i \ol{d}x^i$, put
$$\tau (u,v,x^2,\ldots,x^{n-1})=\sum_{i=2}^{n-1} x^i\left(\int^1_0 \beta_i(u,v,\s x^2,\dots, \s x^{n-1})d\s\right)+ h(u,v),$$
for arbitrary  $h$. Then $\ol{d}\beta=0$ (i.e., $\fr {\p \beta_i}
{\p x^j}-\fr {\p \beta_j} {\p x^i}=0$), yields $\partial_i \tau
=\beta_i$, as required. $\square$

\subsubsection{The covariant derivative $\oln$ and its curvature tensor $\ol{R}$ on $\M$}

\begin{defn}\label{covder} The
{\em covariant derivative} $\oln$ on $\M$ associated to a
Brinkmann decomposition $\{u,v\}$ is the map
\begin{center}
\begin{tabular}{cccc}
 $\oln:$&$\G(T \M)\times\G(T\M)$&$\longrightarrow$&$\G(T \M)$\\
&$(X,Y)$&\textrightarrow& $\oln_X Y$
\end{tabular}
\end{center}
defined at each point $p\in M$ by
$$(\oln_X Y)_p := (\oln^{\bar g}_{\breve X} \breve Y)_{p}$$ 
where $\breve X$ and $\breve Y$ are the restrictions of $X, Y$ to
the leaf $M_{(u(p),v(p))}$ and $\oln^{\bar g}$ the Levi-Civita connection of the first fundamental form $\bar g$ on that leaf. 
\end{defn}

\begin{rema}\label{Rem2}{\em a) Clearly, $\oln$ satisfies the formal properties of a symmetric covariant
derivative, as well as $\oln\ol{g}=0$.

b) The  {\em Christoffel symbols}  $\ol{\G}^i_{jk}$ of $\oln$ in
the basis $\{\p_i\}$  are defined by the relation
$\oln_{\p_j}\p_i=\ol{\G}^k_{ij}\p_k$. Trivially,
they are smooth and invariant by the flow of $\partial_v$.

c) The covariant derivative on $\M$ can be extended to
sections in $T^r_s \M$. For the  coordinate basis
$\{\p_i\}$ and any $T\in\G(T^r_s \M)$ one has:
\beq\nonumber
\oln_m T^{i_1\ldots i_r}_{j_1\ldots j_s}=\p_m(T^{i_1\ldots
i_r}_{j_1\ldots j_s})+\sum_{a=1}^r\ol{\G}^{i_a}_{k m}T^{i_1\ldots
i_{(a-1)} k i_{(a+1)}\ldots i_r}_{j_1\ldots
j_s}-\sum_{b=1}^s\ol{\G}^k_{j_b m}T^{i_1\ldots i_r}_{j_1\ldots
j_{(b-1) }k j_{(b+1)} \ldots j_s}\eeq

d) For $X,Y\in\G(TM)$, $\ol{\n_X Y}\neq \oln_{\ol{X}} \ol{Y}$ in
general (for example, if $X=E_0$ then
$\ol{X}=0$). 
}\end{rema}

The covariant derivative $\oln$ yields a natural {\em curvature
tensor ${\ol{\mathcal{R}}}$ of the foliation $\M$} defined formally as the
usual  curvature of $\oln$:
\begin{equation}\label{olRR}\ol{\mathcal{R}}(X,Y)Z=(\oln_X\oln_Y-\oln_Y\oln_X-\oln_{[X,Y]})Z\in\G(T\M), ~\forall X,Y,Z\in\G(T\M) \end{equation}
and also its derived  {\em Ricci tensor}
$\ol{\mathcal{R}{\hbox{ic}}}$ and {\em scalar curvature $\ol{\mathcal{S}}$ of} $\M$. All of them
satisfy the standard symmetries corresponding to a curvature
tensor.

\begin{defn}\label{foli}
Fixed a Brinkmann decomposition $\{u,v\}$. The foliation $\M$  is
called {\em $u$-Einstein} if $\ol{\mathcal{R}{\hbox{ic}}}=\mu
\ol{g}$ for some function $\mu$ such that $d\mu\wedge du=0$. In particular, when
$\mu$ is constant, $\M$ is  called {\em Einstein}, and when
$\mu\equiv 0$ we say that $\M$ is  {\em Ricci-flat}.

In the case $\ol{\mathcal{R}}= 0$ (respectively $\oln~ \ol{\mathcal{R}}= 0$),  the foliation
$\M$ is  said to be {\em flat} (respectively {\em locally symmetric}).
\end{defn}

Some simple properties follow immediately.
\begin{prop}\label{redfol} Let $(M,g)$ be a Brinkmann space with a fixed Brinkmann decomposition $\{u,v\}$.
\ben
\item If 
$\oln^r\ol{\mathcal{R}}=0$ for some $r> 1$, then $\M$ is
a locally symmetric foliation.
\item If $\M$ is locally symmetric and Ricci-flat, then it is flat.
\item If  $\M$ is flat, the Brinkmann decomposition admits a
chart $\{u,v,y^i\}$ such that the metric $g$ becomes:
$$g=-2du(dv+H du+W_i dy^i)+\delta_{ij}dy^i dy^j.$$
\een
\end{prop}

\ni {\em Proof}. (1) For each $(u,v)$, the Riemannian result which states that $\oln^r\ol{\mathcal{R}}=0$ implies local symmetry can be applied to each leaf $\ol{M}_{(u,v)}$.

(2) Apply Proposition \ref{redlocsym}(3) to each leaf of $\M$.

(3) We sketch a procedure to obtain the required new
coordinates $y^{i}=y^{i}(u,x^j)$, for the sake of completeness.
Let $\gamma$ be an integral curve of $\partial_u$ along the
hypersurface $\Omega=\{v=0\}$, and put for small $U \subset \R^{n-2}$, $0\in U$:
$$
\varphi: ]-\epsilon,\epsilon[ \times U \rightarrow
M , \qquad \qquad (u,w^i)\mapsto
\ol{\exp}_{\gamma(u)}(w^i\partial_i|_{\gamma(u)}),$$ where
$\ol{\exp}_{\gamma(u)}$ denotes the exponential at $\gamma(u)$ on
the leaf $\ol{M}_{(\gamma(u),0)}$ for the (flat)
metric
$\ol{g}$.
The classical Cartan theorem shows that each map $\varphi_u
:=\varphi(u,\cdot)$ is an affine transformation from $U$ to the
leaf. So, regarding $\{ w^i\}$ as coordinates in
$\R^{n-2}$,
$$\varphi_u^*\ol{g} = h_{ij}(u) dw^i dw^j$$ for some {constants} $h_{ij}(u)$ on each leaf
 (which define an Euclidean metric and depend smoothly
on $u$). Now, $\varphi$ allows to consider $\{u, w^i\}$ as
coordinates in $\Omega$, with each coordinate vector field
$\partial/\partial w^i$ parallel on each leaf of $\M_{(u,0)}$.
 By using the Gram-Schmidt procedure,
 an orthonormal basis 
 $\{V_{j}=B_{j}^{i}\frac{\partial}{\partial
 w^i}\}$ in $T\M$ is obtained. Indeed, by construction the transition matrix $(B_j^i)$
 depends smoothly only on $u$ and, thus, by Remark \ref{rlata} (3) $\oln
V_{j}=\ol{d}(B_{j}^{i}) \frac{\partial}{\partial
 w^i}=0$. Therefore, Proposition \ref{df1} (4) implies that $(V_j)^\flat= dy_\Omega^j$
 for some functions $\{y_\Omega^j(u,w^j)\}$ on some open subset of $\Omega$. The required functions $\{y^j\}$ are obtained by extending $\{y_\Omega^j(u,w^j)\}$ to
 a neighborhood of $M$ in a $v$-invariant way according to Remark \ref{rlata} (2).
  $\square$

\subsubsection{The $\dot{}$ (dot) derivative }
We introduce the following simple derivative:
\begin{defn} {\em The dot derivative} $\dot{T}\in \G({T}^r_s \M)$ of a tensor field
$T\in \G({T}^r_s \M)$ is defined as $\dot{T}={(\ol{\mc{L}_{\p_u}\ri T})}.$
\end{defn}
The components of $\dot{T}$ in the coordinate
basis $\{\p_i\}$ are $\dot{T}^{i_1\ldots i_r}_{j_1 \ldots
j_s}=\p_u ({T}^{i_1\ldots i_r}_{j_1 \ldots j_s}).$ Indeed, the first usage of the dot
for functions was the definition \eqref{e0}. 

\subsection{Some geometrical objects and the operator $D_0$}\label{s101}
Our aim in this section is to obtain explicit expressions for the curvature and its derivatives
adapted to a Brinkmann chart. We start by computing the curvature
two-forms for a \nr . Then, we introduce a new operator $D_0$ which, together
with those introduced in the last section, will allow us to simplify calculations and to provide manageable formulae for the derivatives of the curvature. This will be very helpful in order to solve the equations of 2nd-symmetry.

\subsubsection{The connection $1$-forms and the curvature
tensor}\label{sconn1form}
We compute the connection and curvature forms by
using Cartan techniques, according to the conventions in \cite{NAK}\cite{CI}.
In a standard manner, we define the
connection one-forms $\ol{\o}^i_j \in \Lambda^1(\M )$ and the corresponding
curvature two-forms $\ol{\O}^i_j \in \Lambda^2(\M )$
of the foliation $\M$ for the given coordinate
basis $\{\p_i\}$ as:
$$\oln_X \p_i=\ol{\omega}^j_i(X) \p_j, \quad
\ol{\mathcal{R}}(X,Y){\p_{j}}= \ol{\O}^i_j(X,Y)\p_i, \quad \forall X,Y\in\G(T\M).$$
Observe that $\ol{\o}^i_j=\ol{\G}~^i_{jk}\ol{d}x^k$, where
$\ol{\G}~^i_{jk}$ are the Christoffel symbols of $\oln$ as introduced in Remark
\ref{Rem2}. Then, a simple computation
 shows that the  first and second Cartan's
equations  for $\M$ still hold:
\beq \nonumber 0&=&\ol{\o}^i_j\wedge\ol{d}x^j \\\nonumber
\ol{\O}^i_j&=&\ol{d}~\ol{\o}^i_j+\ol{\o}^i_k\we\ol{\o}^k_j \, .
\eeq

In order to compute the connection and curvature of the Brinkmann spaces we introduce two tensor fields on $\M$ in terms of the
adapted differential $\ol{d}$ and the dot derivative, which will
be especially relevant in the computations.
We define $h\in \Lambda^1 ( \M )$ and $t\in T^0_{2}(\M)$ as
\beq \label{h} h&=&\ol{d}H-{\dot{W}},\\
\label{t} t&=&-\frac{1}{2}\left(\dot{\ol{g}}+\ol{d}{W}\right).
\eeq
Observe that the symmetric and skew-symmetric parts of $t$ are precisely $-\dot{\ol{g}}/2$ and $-\ol{d}{W}/2$, respectively.
We emphasize that $h$ and $t$ depend on the Brinkmann chart.
Equipped with these sections, one can check that the non-vanishing connection one-forms associated to the a \nr
$\{E_\a\}$ of any Brinkmann chart are
\beq\label{W1}
 \o^1_i&=&\o^j_0 g_{ij}=h_i \t^0 - t_{ij}\t^j,
\\
\label{W2}
 \o^i_j&=&-g^{ik} t_{kj}\t^0+\ri {\ol{\o}}^i_j.
\eeq
where $h_i=(H_{,i}-\dot{W}_i)$ and $t_{ij}=\fr 1 2 (-\dot{g}_{ij}+W_{i,j}-W_{j,i})$ are the components of the tensors defined in \eqref{h} and \eqref{t}. It should be noted that $\ol{\o^i_j}=\ol{\o}^i_j$, but $\ri {\ol{\o}}^i_j\neq \o^i_j$.
Then, the non-vanishing curvature two-forms associated to the \nr $\{E_\a\}$ of any Brinkmann chart read
\beq
\label{new1}
 \O^1_i&=&-(\ol{\nabla}_j h_i+\dot{t}_{ij}+t^k~_{i}t_{kj})\t^0\we\t^j+\frac{1}{2}(\ol{\nabla}_k t_{ij}-\ol{\nabla}_j t_{ik})\t^j\we\t^k,
\\
\label{new2}
 \O^i_j&=&(\oln{_k}t^i~_j+\dot{\ol{\G}}~^i_{jk})\t^0\we\t^k+\overline{\O}^i_j. 
\eeq
so that again $\ol{\O^i_j}=\ol{\O}^i_j$, but $\ri {\ol{\O}}^i_j\neq \O^i_j$.

From (\ref{W1}-\ref{W2}) it is immediate to obtain the components $\g^\a_{\b\la}$ of the connection one-forms $\o^\a_\b$, defined by $\o^\a_\b=\g^\a_{\b\la}\t^\la$. These will have to be used later. Recall also the identity $\n_{E_\b}E_\a=\g^\la_{\a\b}E_\la$. This together with (\ref{W1}) provides the simple formula
\begin{equation}\label{new3}
h(X)=g(\ri X,\n_{E_0} E_0), \,\,\,\,  \forall X \in\G(T\M).
\end{equation}
Similarly, the non-vanishing components $R^\a{}_{\b
\gamma \delta}$ of the curvature tensor $R$, defined by $\Omega^\a_{\b}=\frac{1}{2}R^\a{}_{\b\gamma \delta}\theta^\gamma \wedge \theta^\delta$ can be read off from (\ref{new1}-\ref{new2}):
\beq
\label{Aij}
 R^1~_{i0j}&=&-(\ol{\nabla}_j h_i+\dot{t}_{ij}+t^k~_{i}t_{kj}),
\\
\nonumber R^1~_{ijk}&=& \ol{\nabla}_k t_{ij}-\ol{\nabla}_j t_{ik},
\\
\nonumber  R^i~_{jkl}&=&\ol{\mathcal{R}}^i~_{jkl},
\\\label{R4}
 R^i~_{j0k}&(=&-g^{ri}R^1~_{krj})=\oln{_k}t^i~_j+\dot{\ol{\G}}{}^i{}_{jk}
\eeq
where, in the last expression, the dot acts on the functions $\ol{\G}{}^i{}_{jk}$.
From here, the non-vanishing components of the Ricci tensor
can be easily computed (note that $R^{1}~_{\a\b\mu}=-R_{0\a\b\mu}$):
\begin{eqnarray}\nonumber R_{00}&=&R^i~_{0i0}=-g^{ij}R^1~_{j0i}=\oln_i h^i+g^{ij}\dot{t}_{ji}+t^{ki}t_{ki},\\
 \nonumber R_{0i}&=&\oln_i t^j~_j-\oln_j t^j~_i,\\
\nonumber R_{ij}&=&\ol{\mathcal{R}}_{ij},
\end{eqnarray}
and the scalar curvature of $M$ turns out to be equal to that of $\M$
$$S=\ol{\mathcal{S}}.$$
Therefore, we have proven the following result.
\begin{prop}\label{R=R}
The curvature tensor $\ol{\mathcal{R}}$ of the foliation $\M$ as defined in \eqref{olRR}, and its associated Ricci tensor $\ol{\mathcal{R}\text{ic}}$ and scalar curvature $\ol{\mathcal{S}}$ satisfy that
$$\ol{\mathcal{R}}=\ol{R};~~~~\ol{\mathcal{R}\text{ic}}=\ol{\text{Ric}};~~~~\ol{\mathcal{S}}=\ol{S},$$
where $\ol{R}$, $\ol{\text{Ric}}$ and $\ol{S}$ are the projections by the homomorphism \, $\bar{}$ \, of the curvature tensor $R$, the Ricci tensor $\text{Ric}$ and the scalar curvature $S$ of the Brinkmann space $(M,g)$, respectively.
\end{prop}
Using this result, from now on there will be no need to distinguish between these objects and we will use only the notation $\ol{R}$, $\ol{\text{Ric}}$ and $\ol{S}$ for them.

\begin{rema}{\rm About the consistency of the tensor equations on $M$:

(1) If, say, the first index of $t$ and $\ri t$
is raised, then the same expression $t^i~_{j}$ is obtained, as the
isomorphisms $\flat$ and $\sharp$ commute with the homomorphism
$T\mapsto \ri T$ ($t^i~_j\equiv g^{i\a}t_{\a j}=g^{ir}t_{rj}$).

(2) Observe that $\ol{R}=\ol{\mathcal{R}}$ is related to $R$ by means of the following formula:
$$\ol{R}(X,Y)Z=\ol{R(\ri X,\ri Y)\ri Z} , \quad\quad \forall X,Y,Z\in\G(T\M).$$
}\end{rema}

\subsubsection{The $D_0$ derivation}
\begin{defn}
The $D_0$ operator associated
to any Brinkmann chart is defined as:
\begin{center}
\begin{tabular}{cccc}
 $D_0:$&$\G(T^r_s \M)$&$\longrightarrow$&$\G(T^r_s \M)$\\
&$T$&\textrightarrow& $D_0 {T}=\ol{(\n_{E_0} \ri T)}$
\end{tabular}
\end{center}
Consequently, ${T}\in\G(T^r_s\M)$ is said to be {\em $D_0$-parallel} if $D_0 {T}=0$.
\end{defn}
$D_0$ measures the variation, projected to $T^r_s\M$, of tensor fields
$T\in\G(T^r_s \M)$ under displacements along the direction $E_0$,
which is transverse to the leaves of $\M$.
\begin{prop}
 $D_0$ satisfies the formal properties of a tensor derivation on $\M$:
\ben \item[i)]$\R$-linearity: $D_0(a A+b B)=a D_0 A+ b D_0 B$,
$\forall a,b \in \R$, $\forall A,B \in\G(T^r_s \M).$ \item[ii)]
Leibniz rule: $D_0(A\otimes B)=(D_0 A)\otimes B+A\otimes(D_0 B)$,
$\forall A,B \in\G(T^r_s \M).$ \item[iii)] Commutativity  with
contractions: $D_0 (C^i_j(A))=C^i_j(D_0 A)$, $\forall A\in\G(T^r_s
\M)$, where $C^i_j$ denotes the contraction of the $i^{th}$ contravariant slot with the $j^{th}$ covariant one. \een
\end{prop}
\begin{prop} For $X\in\G(T\M)$, the decomposition of $\n_{E_0}\ri
X$ in $T\U\oplus \K$ is
$$\n_{E_0}\ri X= h(X) E_1 + \ob{D_0 X}^\circ. $$
\end{prop}
{\em Proof.} Putting $\n_{E_0}\ri X=X_1+X_2$, with $X_1\in \G(T
\U)$ and $X_2\in \G(\K)$, the definition of $D_{0}$ gives immediately $D_0 X=\ol{X_2}$. For $X_1$, using that $g(\ri X,E_0)=g(\ri X,E_1)=0$, that $E_1$ is parallel and formula (\ref{new3})
we get:
$$X_1=-g(\n_{E_0}\ri X,E_1)
E_0-g(\n_{E_0}\ri X,E_0) E_1 =g(\ri X,\n_{E_0} E_1) E_0 + g(\ri
X,\n_{E_0} E_0) E_1 =h(X) E_{1} \, .
\square$$

\begin{prop}
$\ol{g}$ is $D_0$-parallel, that is to say:
\begin{equation}\label{d0g}D_0\ol{g}=0.\end{equation}
\end{prop}

\ni {\it Proof}. By definition, $D_0\ol{g}\equiv\ol{\n_{E_0} \ri
{\ol{g}}}$ and $g=- \t^0\otimes \t^1-\t^1\otimes\t^0+\ri
{\ol{g}}$.
Since $\t^0$ is  parallel and  $\n_{E_0} {g}=0$,
$$-\t^0 \otimes \n_{E_0} \t^1-\n_{E_0} \t^1\otimes \t^0+\n_{E_0}\ri {\ol{g}}=0$$

Therefore, as $\ol{\left(\t^0 \otimes \n_{E_0} \t^1+\n_{E_0}
\t^1\otimes \t^0\right)}=0$, we have that $\ol{\n_{E_0}\ri {\ol{g}}}=0$.
$\square$

\begin{rema}
 {\rm It is not difficult to show that, for arbitrary $\omega\in\G(T^*\M)$, the $D_{0}$-derivative acts such that
$$
(D_{0}\omega)(X) =E_{0}(\omega)(X)+t(\omega^\sharp ,X), \,\,\,\, \forall X\in \G(T\M)
$$
while, with the same notation
$$
\omega (D_{0}X) =\omega (E_{0}(X))-t(\omega^\sharp , X) \, .
$$
These formulas can then be extended to arbitrary sections $T\in\G(T^r_s \M)$, and can be expressed in the local basis $\{\p_i\}$ of $\G(T\M)$ by means of:
\begin{eqnarray*}
(D_0 X)^i&\equiv& D_0 X^i =
(\p_u-H\p_v)(X^i)-t^i~_j X^j,\\
(D_0 \o)_i&\equiv& D_0 \o_i=
(\p_u-H\p_v)(\o_i)+t^j~_i \o_j.
\end{eqnarray*}
Its generalization to any section $T\in\G(T^r_s \M)$ is
$$(D_0 {T})^{i_1\ldots i_r}_{j_1 \ldots
j_s}\equiv D_0 {T}^{i_1\ldots i_r}_{j_1 \ldots
j_s}=(\p_u-H\p_v)({T}^{i_1\ldots i_r}_{j_1 \ldots
j_s})-\sum_{a=1}^r t^{i_a}~_k T^{i_1\ldots i_{(a-1)} k i_{(a+1)}
\ldots i_r}_{j_1 \ldots j_s}+\sum_{b=1}^st^k~_{j_b} T^{i_1\ldots
i_r}_{j_1\ldots  j_{(b-1)} k j_{(b+1)} \ldots j_s}.$$
For $v$-invariant tensor fields the expression $D_0T$ simplifies ($D_0 {T}=\ol{\n_{\p u} \ri T}$) such that
\beq\label{D0T}(D_0 {T})^{i_1\ldots i_r}_{j_1 \ldots
j_s}=\dot{T}^{i_1\ldots i_r}_{j_1 \ldots j_s}-\sum_{a=1}^r
t^{i_a}~_k T^{i_1\ldots i_{(a-1)} k i_{(a+1)} \ldots i_r}_{j_1
\ldots j_s}+\sum_{b=1}^st^k~_{j_b} T^{i_1\ldots i_r}_{j_1\ldots
j_{(b-1)} k j_{(b+1)} \ldots j_s}.\eeq
}
\end{rema}
Next, we collect some elementary properties of the
$D_0$-parallel transport to be used later. First, let $\eta$ be
an integral curve of $E_0$ and take a vector field
$X_\eta$
along $\eta$ which is everywhere tangent to $\M$, i.e.,
${X_\eta}\in\G(T_{\eta}\M)$. The derivative
$D_0(X_\eta)$ of $X_\eta$, as well as its $D_0$-parallelism, makes an obvious sense.
\begin{lema}
Let $p\in M$, $\v v \in T_p \M$ and $\eta$
 the integral curve of $E_0$ with $\eta(0)=p=(\up,\vp,\xip)$. Then, there exists a unique vector field
 $X_\eta$ obtained as the $D_0$-parallel transport along $\eta$ such that ${X}_{\eta(0)}=\v v$.
\end{lema}
A standard reasoning leads to:
\begin{prop}\label{doort}
The map which sends each $\v v\in T_p\M$ to its unique $D_0$-parallel
transport $X_\eta(\tau)\in T_{\eta(\tau)}\M$ along $\eta$ is a linear
isometry from $T_p\M$ to $T_{\eta(\tau)}\M$.
\end{prop}
{\it Proof.} The result follows from \eqref{d0g}, since
$(\ol{g}(X,Y))|_{\eta}$ depends only on the parameter $u$ of $\eta$. $\square$

Our last result yields an extension to all the leaves of any vector field on one leaf.
\begin{prop}\label{exten}
Let $X_{\ol{M}}$ be a vector field on a leaf $\ol{M}_{(u_0,v_0)}$ of $\M$.
Then, there exists a unique $v$-invariant and $D_0$-parallel vector field $X_\M\in\G(T\M)$ which extends $X_{\ol{M}}$.
\end{prop}

{\em Proof}. Consider the
hypersurface $\Omega_{v_0}=\{v=v_0\}$ and extend $X_{\ol{M}}$
 to a vector field $X_{\Omega_{v_0}}$ on
$\Omega_{v_0}$ by taking each integral curve $\eta$ of $E_0$ which
starts at some $p\in \ol{M}_{(u_0,v_0)}$, and defining
$X_{\Omega_{v_0}}\circ \eta$ as the $D_0$-parallel transport of
$X_{p}$. Then, extend $X_{\Omega_{v_0}}$ to a $v$-invariant vector
field according to Remark \ref{rlata}(2). $\square$

\subsubsection{Derivatives of the curvature tensor $R$}\label{expR}
 As $E_1$ is parallel, the curvature tensor will be determined on each
Brinkmann chart by its value on quadruples of vectors tangent to $\M$,
plus some extra tensors which take care of the remaining components (partly) along the $\theta^1$ or $E_0$ directions.
We start by defining two such tensor fields
 on $\M$.
\begin{defn}\label{AB} For any Brinkmann chart and its associated \nr $\{E_\alpha\}$
we  define $A\in \G(T^0_2\M)$ and $B\in \G(T^0_3\M)$ as $A:= \theta^1(R(E_0,\cdot ))$ and $B:=\theta^1 (R)$, that is to say
$$A(X,Y)={\theta^1(R(E_0,\ri Y)\ri X)}, \hspace{1cm}  B(X,Y,Z)={\theta^1(R(\ri Y,\ri Z)\ri X)}\hspace{1cm} \forall X,Y,Z\in\G(T\M).$$
\end{defn}
From the symmetries of the curvature tensor it is obvious that $A$ is symmetric
$$
A(X,Y) = A(Y,X) \hspace{1cm} \forall X,Y,\in\G(T\M)
$$
and that $B$ is skew-symmetric in its last two slots and it  satisfies a cyclic identity
$$
B(X,Y,Z)=-B(X,Z,Y), \hspace{.3cm} B(X,Y,Z)+B(Y,Z,X)+B(Z,X,Y)=0, \hspace{.3cm} \forall X,Y,Z\in\G(T\M).
$$

\begin{rema}\label{rema222} {\rm In the given basis $\{\p_i\}$, $A_{ij}=A(\p_i,\p_j)=R^1~_{i0j}$ and $B_{ijk}=B(\p_i,\p_j,\p_k)=R^1~_{ijk}$. Then, the previous properties can be expressed using this notation as $A_{ij}=A_{(ij)}$, $B_{ijk}=B_{i[jk]}$, and $B_{[ijk]}=0$. In what follows, and for the sake of brevity, we will resort to using index notation in many cases, which is sufficient to illustrate these properties and reveals itself as very helpful in the required complicated calculations for 2nd-symmetry. As a starting example, note that we additionally have for instance
$$B_{ij}~^k (=g^{kr} B_{ijr})=R^k~_{j0i}.$$ }\end{rema}

A direct computation of $\n R$, for instance in the basis $\{E_{\a}\}$, provides the following formulae:
\begin{eqnarray*}
\ol{\theta^1(\n_{E_0} R(E_0,\cdot))}&=&D_0 A+2{\mathcal S} \left[C^1_3 (h^\sharp\otimes B)\right]\\
\ol{\theta^1(\n R(E_0,\cdot))}&=&\ol\n A-2{\mathcal S} \left[C_{15} (t\otimes B)\right]\\
\ol{\theta^1(\n_{E_0} R)}&=&D_0B+h(\ol R)\\
\ol{\theta^1(\n R)}&=&\ol\n B -C^1_1(t\otimes \ol R)\\
\ol{\n_{E_0} R} &=& D_0 \ol R\\
\ol{\n R} &=& \ol\n~\ol R
\end{eqnarray*}
where ${\mathcal S}[T]$ gives the symmetric part of any covariant section $T\in \G(T^0_s\M)$ in its last two slots, and $C_{ij}$ denotes the contraction of the $i^{th}$ and $j^{th}$ covariant indices (via the metric $\bar{g}^{-1}$).


We give names to the lefthand sides of these relations (except for the last one), thereby defining five tensor fields on $\M$ which will allow for simpler expressions when computating $\n^2 R$.
\begin{defn}\label{nRs} For any Brinkmann chart and its associated \nr $\{E_\alpha\}$ 
we  define $\wt{A}\in \G(T^0_2\M)$, $\wh{A},\wt{B}\in \G(T^0_3\M)$, $\wh{B}\in\G(T^0_4\M)$ and $\wt{R}\in\G(T^1_3\M)$, for all $X,Y,Z,V\in\G(T\M)$ by:
\begin{eqnarray*}
 \wt{A}(X,Y)&=&{\theta^1\left((\nabla_{E_0}R)(E_0,\ri Y)\ri X\right)},\\
\wh{A}(X,Y,Z)&=&{\theta^1\left((\nabla_{\ri X}R)(E_0,\ri Z)\ri Y\right)},\\
\wt{B}(X,Y,Z)&=&{\theta^1\left((\nabla_{E_0}R)(\ri Y,\ri Z)\ri X\right)}, \\
\wh{B}(X,Y,Z,V)&=&{\theta^1\left((\nabla_{\ri X}R)(\ri Z,\ri V)\ri Y\right)},\\
\wt{R}(X,Y)Z&=&\ol{\nabla_{E_0} R(\ri X,\ri Y)\ri Z}.
\end{eqnarray*}
\end{defn}
Observe that, therefore, one has $\tilde R =D_0 \ol R$, $\tilde B =D_0 B+h(\ol R)$, etcetera in agreement with the previous formulas.

\begin{rema} {\rm In the given basis $\{\p_i\}$, $\wt{A}_{ij} = \wt{A}(\p_i,\p_j)=\nabla_0 R^1~_{i0j}$, $\wh{A}_{sij}  =  \wh{A}(\p_s,\p_i,\p_j)=  \nabla_s R^1~_{i0j}$, $\wt{B}_{ijk}  = \wt{B}(\p_i,\p_j,\p_k)= \nabla_0 R^1~_{ijk}$, $ \wh{B}_{sijk}  = \wh{B}(\p_s,\p_i,\p_j,\p_k)=  \nabla_s R^1~_{ijk}$ and $\wt{R}^i~_{jkl} = \ol{d}x^i\left(\wt{R}(\p_k,\p_l)\p_j\right) = \nabla_0 R^i~_{jkl}.$ Besides, from the symmetries of the curvature tensor it is obvious that
\begin{enumerate}
\item $\wt{A}$ and $\wh{A}$ are symmetric in the last two indices: $\wt{A}_{ij}=\wt{A}_{(ij)};~~\wh{A}_{sij}=\wh{A}_{s(ij)}.$
\item $\wt{B},\wh{B},\wt{R}$ are skew-symmetric in their last two indices: $\wt{B}_{ijk}=\wt{B}_{i[jk]};\wh{B}_{sijk}=\wh{B}_{si[jk]};\wt{R}_{sijk}=\wt{R}_{si[jk]}$.
\item  $\wt{B},\wh{B},\wt{R}$ satisfy a cyclic identity: $\wt{R}_{i[jkl]}=0;\wt{B}_{[ijk]}=0;\wh{B}_{s[ijk]}=0$.
\item $\wt{R}$ also satisfies that $\wt{R}_{ijkl}=-\wt{R}_{jikl}$ and $\wt{R}_{ijkl}=\wt{R}_{klij} $, so that it has all the symmetries of a Riemann tensor.
\end{enumerate} }
\end{rema}
One can prove that, in addition, $\wt{B}_{ij}~^k=\nabla_0 R^k~_{j0i}$ and
$\wh{B}_{sij}~^k=\nabla_s R^k~_{j0i}.$ We also point out the following basic relations:
\beq\label{b1}
\wt{R}_{ijkl}&=&-2\wh{B}_{[ij]kl}\\
\label{b2}
 \wt{B}_{kij}&=&2\wh{A}_{[ij]k}\\
\nonumber
\wh{B}_{[i|l|jk]}&=&0 \hspace{17mm} ( \Longrightarrow \wh{B}_{[i|l|j]k}=-\frac{1}{2}\wh{B}_{klij}\, )
\eeq
They follow by direct application of the second Bianchi identity $\nabla_{[\a}R_{\b\s]\r\mu}=0$ in the \nr $\{E_\a\}$: the first two relations follow by taking $\{\a,\b,\s\}=\{0,i,j\}$ and
the last one with $\{\a,\b,\s\}=\{i,j,k\}$.

Using all the above, another direct computation of $\n\n R$
leads to the following formulae:
\begin{subequations}\label{R}
 \begin{align}
 &\ol{\n\n R}=\ol\n\, \ol\n\, \ol R,  & &\ol{\n_{E_0} \n R}=D_0\oln~\ol R,\\
&\ol{\n \n_{E_0}R}={\oln~\wt{R}}+C_{13}(t\otimes\oln~\ol{R}),  & & \ol{\n_{E_0} \n_{E_0} R}=D_0 \wt{R}-C^1_1 (h^\sharp\otimes\oln\, \ol{R}),\\
&\ol{\t^1(\n\n R)}=\oln \wh{B}-C^1_1(t\otimes\oln~\ol{R}), &  & \ol{\t^1(\n_{E_0}\n_s R)}=D_0 \wh{B}+C^1_1(h\otimes\oln ~\ol{R}),\\
& \ol{\t^1(\n \n_{E_0}R)}=\oln \wt{B}-C^1_1 (t\otimes\wt{R})+C_{13}(t\otimes \wh{B}),& &\ol{\t^1(\n_{E_0}\n_{E_0} R)}=D_0 \wt{B}+C^1_1(h\otimes\wt{R})-C^1_1(h^\sharp\otimes\wh{B}),\\
&  \ol{\t^1(\n\n R(E_0,\cdot))}=\oln \wh{A}-2{\mathcal S}[C_{16}(t\otimes\wh{B})],& & \ol{\t^1(\n_{E_0}\n R(E_0,\cdot))}=D_0 \wh{A}+2 {\mathcal S}[C^1_4(h^\sharp\otimes\wh{B})],
\end{align}
\begin{align}
 &\ol{\t^1(\n \n_{E_0} R(E_0,\cdot))}=\oln \wt{A}-2{\mathcal S} [C_{15}(t\otimes\wt{B})]+C_{13}(t\otimes\wh{A}),\\
&\ol{\t^1(\n_{E_0}\n_{E_0}R(E_0,\cdot))}=D_0 \wt{A}+2{\mathcal S} [C^1_3(h^\sharp\otimes\wt{B})]-C^1_1(h^\sharp\otimes\wh{A}).
 \end{align}
\end{subequations}
Recall that, for example, the equations in the first row can be written as:
\begin{align}
\label{lsR} \nabla_m\nabla_s R^i~_{jkl} &= \oln_m\oln_s \ol{R}^i~_{jkl}, &
\nabla_0\nabla_s R^i~_{jkl} &= D_0 \oln_s \ol{R}^i~_{jkl}.
\end{align}


\subsection{Reducibility and a generalized Eisenhart theorem}\label{s7a}
The proof of our main Theorem \ref{lem4} will be carried out by means of some local
decompositions of the Brinkmann space. To that end, we need an appropriate notion of
reducibility and a version of a classical theorem by Eisenhart 
adapted to Brinkmann spaces. 
\begin{defn}\label{decg} A Brinkmann decomposition $\{u,v\}$ of a Brinkmann space $(M,g)$ is {\em spatially reducible}  if there exists a Brinkmann chart
$\{u,v,x^i\}$ around each point and a partition of the indices
$I_1=\{2,\dots, d+1\}, I_2=\{d+2,\dots, n-1\}$ for some $d\in
\{1,\dots , n-3\}$ such that $\ol{g}_{aa'}=0$ and ${\p_{a'}\ol{g}_{ab}}=0$, where the unprimed indices $a,b$ always
belong to the same subset $I_m$ ($m\in\{1,2\}$) and the primed ones $a',b'$
to the other one.

For such a Brinkmann chart, we say that $T\in\G(T^r_s\M)$ is {\em reducible} whenever $T=T^{(1)}+T^{(2)}$, with
$T^{(m)}={T^{(m)}}^{a_1\ldots a_r}_{b_1\ldots b_s}(u,x^c) \p_{a_1}\otimes\ldots\otimes\p_{a_r}\otimes \ol{d}x^{b_1}\otimes\ldots\otimes\ol{d}x^{b_s}$, and $ a_1,\ldots,a_r,b_1,\ldots,b_s,c\in I_m$. We will denote such decomposition as $T=T^{(1)}\oplus T^{(2)}$.
\end{defn}
\begin{rema}\label{*}{\rm (1) Observe that, if a Brinkmann decomposition is spatially reducible, then the metric $\ol{g}$ on $T\M$ is reducible as a tensor field.

(2) Spatial reducibility implies the following property, which provides a more intrinsic expression for some of its consequences. For some Brinkmann
decomposition $\{u,v\}$ there exist two foliations $\M^{(1)},
\M^{(2)}$ of $M$ such that $T\M = T\M^{(1)} \oplus T\M^{(2)}$, and
this sum is orthogonal with respect to $\ol{g}$. In this case, we write $\M= \M^{(1)}\times
\M^{(2)}$ and $\ol{g}=\ol{g}^{(1)}\oplus \ol{g}^{(2)}$, according to the notation in Definition \ref{decg}. Then, the metric $g$ on $M$ can be writen as:
$$g=-2du(dv+H du + \ri W)+\ri {\ol{g}}^{(1)}\oplus\ri {\ol{g}}^{(2)}$$
(recall that $H$ and $ W$ depend on the chosen coordinates $\{x^i\}$). In this context, a tensor field $T$ is reducible if and only if $T=T^{(1)}+ T^{(2)}$ where each $T^{(m)}$ is invariant by the flow of vectors in $\M^{(m')}$ ($\mc{L}_X T^{(m)}=0, \forall X\in\G(T\M^{(m')})$) and it vanishes when applied on any element of $T\M^{(m')}$ and $T^*\M^{(m')}$. In particular, this happens for $\ol{g}$.
The Riemannian manifold $(\ol{M},\ol{g})$ can be written as the product of two
manifolds which will also be denoted, abusing the notation, by
$(\ol{M}^{(1)},\ol{g}^{(1)})$ and $(\ol{M}^{(2)},\ol{g}^{(2)})$,
each $\ol{M}^{(m)}$ generating $\M^{(m)}$ as $\ol{M}$ generated $\M$, see Section \ref{s100}. Thus, the equality
$\ol{g}=\ol{g}^{(1)}\oplus \ol{g}^{(2)}$ may refer either to the
metric decomposition in a leaf $\ol{M}_{(u_0,v_0)}$ 
or in $\M$. The latter depends on $u$ and is $v$-invariant.
Even though this ambiguity is harmless, we will always refer to
decompositions in $\M$ except if otherwise is specified.

(3) Sometimes, a Brinkmann
chart may admit a partition $I_1, \dots , I_s, s\geq 2$ of the
indices $\{2,\dots , n-1\}$ so that 
$\bar g \in \Gamma(T^0_2\M )$ satisfies the properties given in Definition
\ref{decg} for each $i,j\in I_m$, $k',l'\in I_{m'}$ and $m\neq m'$. In this case,  the notation for orthogonal decomposition is naturally extended:
$${\ol{M}}= {\ol{M}}^{(1)}\times \ldots\times
{\ol{M}}^{(s)}, \quad {\M}= {\M}^{(1)}\times \ldots\times
{\M}^{(s)}, \quad  {\ol{g}}= {\ol{g}}^{(1)}\oplus\ldots\oplus
{\ol{g}}^{(s)},$$ and notions as being Einstein or flat are also
extended to each ${\ol{M}}^{(m)}$, ${\M}^{(m)}$ in a trivial way. However, the following caution
must be kept in mind. The given decomposition of $\M$ induces
also an orthogonal decomposition at each leaf $\ol{M}_{(u,v)}$, in
particular at $\ol{M}$. Nevertheless, as $\ol{g}$ is
``$u$-dependent'' such a decomposition may be irreducible in the
sense of the traditional de Rham's theorem for some leaves, but reducible for
other ones. In principle, we will not care about this possible
``spatial irreducibility'' of the metrics $\ol{g}^{(m)}$ on the leaves. Eventually, though,
we will arrive at a decomposition of $\ol{g}$ which will induce an irreducible decomposition of all the leaves, independent of $u$.}
\end{rema}

Let us turn to the Eisenhart theorem. Its classical version \cite{Eis} states:
\begin{teor}\label{Eise1} If a Riemannian manifold $(N,g_{R})$ admits a symmetric two-covariant tensor field $L\in \Gamma(T^0_2 N)$ not proportional to the metric $g_{R}$ such
that $\nabla^{g_{R}} L=0$, then
\begin{itemize}
 \item $g_{R}$ is reducible: $g_{R}=g_{R}^{(1)}\oplus g_{R}^{(2)}\oplus\ldots\oplus g_{R}^{(s)}$ (with each $g_{R}^{(m)}$ not necessarily irreducible).
\item $L=\sum_{m=1}^s \la_m~ g_{R}^{(m)}$ for some constants $\la_m$.
\end{itemize}
\end{teor}
Our aim now is to prove a version of this theorem adapted to the
spatial reducibility of Definition \ref{decg} for Brinkmann decompositions. 
In our generalized version, the reduced metrics $\ol g^{(m)} $ will be
dependent on $u$ but the $\la_m$ will still be constants, 
independent of $u$.
\begin{teor}\label{gei}
Let $(M,g)$ be a Brinkmann space and  fix a Brinkmann chart $\{u,v,x^i\}$. 
 Assume that
there exists a symmetric $v$-invariant,
 $\oln$-parallel and $D_0$-parallel section $\ol{L}\in\Gamma(T^0_{2}\M)$ which is
not proportional to $\ol{g}$. Then, the Brinkmann decomposition $\{u,v\}$ is spatially reducible and $\ol{L}$ is reducible. Furthermore, the decomposition $\{u,v\}$ admits a Brinkmann chart  $\{u,v,y^i\}$ such that:
\ben
 \item ${\ol{g}}= {\ol{g}}^{(1)}\oplus\ldots\oplus
{\ol{g}}^{(s)}$ for some $s\geq 2$.
\item $\ol{L}=\sum_{m=1}^s \la_m~
\ol{g}^{(m)}$ for some constants $\la_m \in \R$.
\een
\end{teor}

{\it Proof}.
Let $p$ be any
point of the chart. We construct an orthonormal basis of eigenvector fields $V_i\in T\M$ defined on the hypersurface $\Omega_{\vp}=\{v=\vp\}$ as we describe now: consider the eigenvalue problem of $\ol{L}_{p}$
with respect to $\ol{g}_p$ on the vector space $T_{p} \M$, i.e.,
\beq\label{eeigen} \ol{L}_{p}(\ldotp,\v v)-\la \ol{g}_{p}(\ldotp,\v
v)=0, \eeq
and take an orthonormal basis $\{\v v_i\}_{i=2}^{n-1}$ of eigenvectors of $\ol{L}_p$ in $T_p\ol{M}_{(\up,\vp)}$. Extend this basis to a normal neighborhood $\ol{U}$ of $p$ in the leaf $\ol{M}_{(\up,\vp)}$ by defining
$V_{i}|_q$ at each $q\in \ol{U}$ as the vector obtained by
 $\oln$-parallelly transporting $\v v_i$ along the unique geodesic
$\gamma_q:[0,1]\rightarrow \ol{U}$ from $p$ to $q$. Clearly, if $\v v_i$ is a $\la$-eigenvector, then $V_{i}|_q$
is an eigenvector of $\ol{L}_{q}$ with the same eigenvalue $\la$, because the one-forms on
$\gamma_q$ defined as $\tau\mapsto
\ol{L}_{\gamma_q(\tau)}(\ldotp,V_{\gamma_q(\tau)})$ and $\tau\mapsto
\la\ol{g}_{\gamma_q(\tau)}(\ldotp,V_{\gamma_q(\tau)})$ are parallel and coincide at $p$ due to
\eqref{eeigen}. Besides, $\{V_i\}$ is an orthonormal basis on $\ol{U}$.  Now, obtain
the sought basis $\{V_i\}$ on $\Omega_{\vp}$ by propagating 
each $V_{i}|_q$ in a $D_0$-parallel manner along the integral curve of $E_0$ at $q\in
\ol{U}$. 
Since $\ol{L}$ is $D_0$-parallel, $\{V_{i}\}$
is still an orthonormal basis of eigenvector fields of $\ol{L}$ and the eigenvalues of $\ol{L}$ are constant on $\Omega_{\vp}$. 

From $v$-invariance, the eigenvalues and the dimension of each eigenspace remain constant all over $M$. Therefore, if we denote any $\la$-eigenvector field as  $V^{(\la)}$, and $\la_1, \dots , \la_s$ are the eigenvalues of $\ol{L}$ with corresponding multiplicities $m_1, \dots
, m_s$, we can reorder the vector fields so that
$$\{V^{(\la_1)}_{1},\ldots\,V^{(\la_1)}_{m_{1}},\ldots\,V^{(\la_s)}_{1},\ldots,V^{(\la_s)}_{m_{s}}\}$$
is an orthonormal basis of $T\M$.

Let $\la$ be one of the eigenvalues and let us prove that the
distribution $S_\la$ generated by its eigenvectors is involutive. Taking the $\oln_{V}$ covariant
derivative of 
$\ol{L}(\cdotp,V_i^{(\la)})=\la
\ol{g}(\cdot,V_i^{(\la)})$ 
 for any ${V} \in \Gamma(T\M)$ and using that $\ol{L}$ is
$\oln$-parallel, we have that $\nabla_{V} V_i^{(\la)}$ lies in
$S_\la$,
and so does $[V_i^{(\la)},V_j^{(\la)}]$. Analogously, the distribution
$S_\la^\perp$ which assigns to each point $p'$ 
the orthogonal of $(S_\la )_{p'}$ in $T_{p'}\M$ is involutive.
Regarding $S_\la^\perp$ as a distribution contained in $T\Omega$,
there are $m_\la +1$
 functionally independent  functions  on $\Omega$ which are
 solutions of the equation $X(f)=0, \forall X\in S_\la^\perp$. The first of
 these functions can be chosen as $u|_{\Omega}$ for all $\la$.
 The other functions
 $y^{i}_\la, i=1,\dots m_\la$, will complete a cooordinate chart
 for $\Omega$ when the construction is repeated for all the
 eigenvalues $\la=\la_j, l=1,\dots , s$. From these coordinates $\{u|_\Omega,y^{i_j}_{\la_j}: i_j=1,\dots
 m_{\la_j}, j=1,\dots , s\}$ on $\Omega$ we can construct a
 chart on $M$ by extending the previous functions in a
 $v$-invariant way according to Remark \ref{rlata} thereby including the
 coordinate $v$.
$\square$

\begin{rema} {\rm
Clearly, under the hypotheses of Theorem \ref{gei}, $\ri {\ol{L}} \in
\Gamma(TM)$ is also diagonalizable (observe that $\partial_u$ and
$\partial_v$ are associated to the eigenvalue 0).
However, the hypotheses on $\ol{L}$ are not enough to ensure that $\ri
{\ol{L}}$ is parallel for $(M,g)$.}
\end{rema}

\section{Proper $2$nd-symmetric Lorentzian manifolds}\label{s6}
We are now equipped with all the elements that will allow us to solve the problem of 2nd-symmetry. The only remaining task is to solve the equations for 2nd-symmetry, given by setting all the expressions in \eqref{R} equal to zero. We are going to do this in several steps.

\subsection{Reduction of the equations}\label{s32}
We start by proving a fundamental simplification of $2$nd-symmetric Brinkmann spaces (that is, of all proper 2nd-symmetric spaces), interesting in its own right.
\begin{teor}\label{redu}
 Let $(M,g)$ be a Brinkmann space with a fixed Brinkmann decomposition $\{u,v\}$. Then, if $(M,g)$ is a $2$nd-symmetric manifold, it follows that:
\ben
\item the foliation $\M$ is locally symmetric, i.e., $\oln~ \ol{R}=0$
\item for any associated Brinkmann chart, the tensor fields $\wh{B}$, $\wt{R}$, $\wh{A}$, $\wt{B}$ in Definition \ref{nRs}  vanish,
\item the scalar curvature $S$ of the manifold is constant. 
\een
\end{teor}
The proof of this theorem will be carried out in two steps, one involving the intrinsic geometry of $\M$ (beginning of Section \ref{first}), and the other involving the integrability equations (Propositions \ref{prop20} and \ref{prop21} in Section \ref{second}). An important consequence of this theorem, to be used later, is:
\begin{coro}\label{Aten}
 Let $(M,g)$ be a 2nd-symmetric Brinkmann space with a fixed Brinkmann decomposition $\{u,v\}$. Then, the section $\wt{A}$  is a tensor field on $\M$, independent of the chosen Brinkmann chart.
Moreover, $(M,g)$ is \tsns if and only if $\wt{A}\neq 0$.
\end{coro}
{\it Proof.} From Theorem \ref{redu}, the only non-zero components of $\n R$ in any \nr associated to a Brinkmann chart $\{u,v,x^i\}$ are $\nabla_0 R^1~_{i0j}(=\wt{A}_{ij})$. Then, it is straightforward to check that $$\wt{A}_{i'j'}=\fr {\p x^{i}} {\p x^{i'}}\fr {\p x^{j}} {\p x^{j'}}\wt{A}_{ij}$$ under a change of the \nr associated to a transformation of the type \eqref{chcoo}. In particular, $\wt{A}$ behaves as a tensor field on $\M$ for the given Brinkmann decomposition $\{u,v\}$.$\square$
\subsubsection{First step: local symmetry of $\M$}\label{first}
To prove (1) of Theorem \ref{redu}, observe that the first in \eqref{lsR} (when set to zero) together with Proposition \ref{redfol} imply the result immediately. At this stage, we can also prove (3) {\em if} we assume (2): given that $S=\ol{S}$ and due to $\oln~\ol{R}=0$ the function $\ol{S}$ depends only on $u$. Thus, assuming $\wt{R}=0$, which is part of  (2), and recalling that $\wt R =D_0 \ol R$ it follows that $0=D_0 \ol{S}$, so that $S=\ol S $ is constant, as required. 

Using \eqref{R} and Theorem \ref{redu}(1), the equations of $2$nd-symmetry thus become

\begin{subequations}\label{Rv2}
 \begin{align}
\label{v21}&\oln_n\wh{B}_{sijk} = 0, & &D_0 \wh{B}_{sijk}   = 0,\\
\label{v22}&\oln_s \wt{R}^i~_{jkl}= 0,&  &D_0 \wt{R}^i~_{jkl}  = 0,\\
\label{v23}&\oln_k \wh{A}_{sij} = 2t^r~_k \wh{B}_{s(ij)r}, &  &D_0 \wh{A}_{sij}  = -2h^r\wh{B}_{s(ij)r},\\
\label{v24}&\oln_s \wt{B}_{ijk} = t^r~_s (\wt{R}_{rijk}-\wh{B}_{rijk}), & &D_0 \wt{B}_{ijk}= -h^r (\wt{R}_{rijk}-\wh{B}_{rijk}),\\
\label{se2}&\oln_k \wt{A}_{ij} = t^r~_k(2\wt{B}_{(ij)r}-\wh{A}_{rij}),\\
\label{se1} &D_0 \wt{A}_{ij} = -h^r(2\wt{B}_{(ij)r}-\wh{A}_{rij}),\\
\label{se11} & \oln_s \ol{R}^i~_{jkl} =0,
 \end{align}
where, for completeness, we include the expressions for all these objects in this notation:
\begin{align}
\label{sss11}
&\wt{R}^i~_{jkl}  = D_0 \ol{R}^i~_{jkl}\\
\label{q4}
&\wh{B}_{sijk}= \oln_s B_{ijk} -t_{rs}\ol{R}^r~_{ijk},\\
\label{qe3}& \wt{B}_{ijk} = D_0 B_{ijk} + h_r \ol{R}^r~_{ijk},\\
\label{qe2}&\wh{A}_{sij}  = \oln_s A_{ij}-2 t^k~_sB_{(ij)k},\\
\label{qe1}&\wt{A}_{ij} = D_0 A_{ij} + 2 h^k B_{(ij)k}\\
\label{BBB}& B_{ijk} = \oln_k t_{ij}-\oln_j t_{ik},\\
\label{AAA}& A_{ij} =  t_{ir} t^r~_j-\oln_j h_i-D_0 t_{ij},\\
\label{RRR}&\ol{R}^i~_{jkl}=R^i~_{jkl}.
\end{align}
\end{subequations}

\subsubsection{Auxiliary algebraic results} \label{s31}
As an interlude, before starting with the second step we
prove a couple of technical results about tensors with particular properties on a real vector space $\mathcal{V}$ of finite dimension. Tensors will be denoted in abstract index form.
\begin{prop}\label{lm2}
 Let ${\mathcal{V}}$ be an $l$-dimensional vector space with a positive definite inner product and $T_{ijk}$ a three-covariant tensor such that:
\ben
\item[(a)]\label{Tt1}  it is skew-symmetric in the last two indices: $ T_{i[jk]}=T_{ijk},$
\item[(b)]\label{T2} it satisfies a cyclic identity: $T_{ijk}+T_{jki}+T_{kij}=0$
\een
If  $T_{(ij)k}$ satisfies \beq \label{TT}T_{(ij)}~^r T_{rnm}=0\eeq then $T_{ijk}=0$.
\end{prop}

{\it Proof}. Observe that, on using (a), (b) can be rewritten as $T_{[ijk]}=0$.
We use arguments inspired in \cite[Lemma 4.1]{SN1}.
Suppose first that there does not exist a non-vanishing vector $v\in \mathcal{V}$ such that $v^r T_{rnm}=0.$ Define the set of vectors $\{Q^r(i,j)\}_{i,j=1}^{l}$ by $Q^r(i,j)=T_{(ij)}~^r$ for each pair $(i,j)$. Then, equation \eqref{TT} can be rewritten as $Q^r(i,j)T_{rnm}=0$, that is, $Q^r(i,j)=0$ for all $i,j$ by our assumption. Consequently, $T_{(ij)k}=0$, i.e., $T_{ijk}=T_{[ij]k}$. Using this fact and (a), $T_{ijk}$ is totally skew-symmetric, and so it vanishes by (b).

Suppose then that there does exist a non-vanishing vector $v\in
\mathcal{V}$ such that $v^r T_{rnm}=0.$ Without loss of
generality, we can assume that $||v||=1$. The orthogonal splitting
of $T_{ijk}$ with respect to $v$
is, then: \beq\label{ort}T_{ijk}=a_{ijk}+b_{ij}v_k-b_{ik}v_j\eeq
where
$$a_{[ijk]}=0,~~a_{i[jk]}=a_{ijk};~~a_{ijk}v^i=a_{ijk}v^j=0,~~b_{ij}=T_{ijk}v^k;~~b_{ij}v^i(=T_{ijk}v^k
v^i)=0,~~b_{ij}v^j=0$$ Since (b) implies $v^i
(T_{ijk}+T_{jki}+T_{kij})=0$, we obtain by \eqref{ort} that
$b_{[jk]}=0$, i.e., $b_{ij}$ is a symmetric two-covariant tensor.
But \eqref{TT} implies $v^i v^mT_{(ij)}~^r T_{rnm}=0$, which in
turn implies $ b_{j}~^r b_{rn}=0$ on using \eqref{ort} once more.
Contracting the indices $j,n$ and using that $b_{ij}$ is
symmetric, we have that $b_{ij} b^{ij}=0$, and as the inner
product is positive definite, that $b_{ij}=0$. Hence,
$T_{ijk}=a_{ijk}$ so that $T_{ijk}$ must be actually totally
orthogonal to $v$. As the tensor $a_{ijk}$ satisfies the
symmetries (a) and (b) of $T_{ijk}$ but in the $(l-1)$-dimensional
space $<v>^\perp$, the result follows by induction. $\square$

An immediate consequence, to be used later, is
\begin{coro}\label{lm1}
 Let ${\mathcal{V}}$ be an $l$-dimensional vector space with a positive definite inner product and $T_{ijkl}$ a four-covariant tensor such that:
\ben
\item[(a)] it is skew-symmetric in the last two indices: $T_{si[jk]}=0,$
\item[(b)] it satisfies the cyclic identity: $T_{ijkl}+T_{iklj}+T_{iljk}=0$
\een
If $T_{s(ij)}~^r T_{lrnm}=0$ holds, then $T_{sijk}=0$.
\end{coro}
{\it Proof}.
Let us define for any $v\in \mathcal{V}$ the tensor $v^i T_{ijkl}$. This tensor satisfies all the conditions in Proposition \ref{lm2}, so  $v^i T_{ijkl}=0$ for all $v\in \mathcal{V}$, that is, $T_{ijkl}=0$. $\square$

\subsubsection{Second step: integrability equations}\label{second}
Obviously, $2$nd-symmetry implies semi-symmetry. As a matter of fact, the equations of semi-symmetry happen to be some of the integrability conditions for the equations of $2$nd-symmetry. To check this, and in order to exploit the integrability condition in full, we introduce the expressions for the commutativity of $\oln$ and $D_0$.

\begin{prop}
 Given a fixed Brinkmann chart $\{u,v,x^i\}$:
\ben
\item[(1)]  for any $T\in \G(T^r_s\M)$, the Ricci identity reads:
\beq\label{Riciden} (\oln_n\oln_s-\oln_s\oln_n) T^{i_1\ldots 
i_r}_{j_1\ldots j_s}=\sum_{b=1}^s \ol{R}^k~_{j_b ns}T^{i_1\ldots
i_r}_{j_1\ldots j_{b-1} k j_{b+1} \ldots j_s}-\sum_{a=1}^r
\ol{R}^{i_a}~_{kns}T^{i_1\ldots i_{a-1} k i_{a+1}\ldots
i_r}_{j_1\ldots j_s} \eeq
\item[(2)] the
commutation of $D_0$ with $\oln$ for any
$F\in \G(T^1_1\M)$ (trivially extendable to arbitrary $T\in
\G(T^r_s\M)$) is:
\begin{equation}
\nonumber(\oln_k D_0-D_0 \oln_k)F^i~_j=(H_{,k}) (\p_v F^i~_j)+F^i~_r B_{kj}~^r-F^r~_j B_{kr}~^i-t^r~_k\oln_r F^i~_j .\end{equation}
\noindent and, when $F$ is $v$-invariant this simplifies to:
\begin{equation}\label{conmF}(\oln_k D_0-D_0 \oln_k)F^i~_j=F^i~_r B_{kj}~^r-F^r~_j
B_{kr}~^i-t^r~_k\oln_r F^i~_j \, .\end{equation}
\een
\end{prop}

{\it Proof}.
(1) follows from the Ricci identities for $\oln$ on each leaf of
$\M$.
To prove
(2), using \eqref{D0T} and Remark \ref{Rem2} we get
\beq\nonumber\oln_k (D_0 F^i~_j)&=&\oln_k(\p_u F^i~_j- H \p_v
F^i~_j-t^i~_r F^r~_j+t^r~_j F^i~_r), \\
\nonumber D_0 (\oln_k F^i~_j)&=&D_0(\p_k
F^i~_j+\ol{\G}^i_{rk}F^r~_j-\ol{\G}^r_{jk}F^i~_r).\eeq The result
follows by expanding these two expressions, and taking into account
$\p_v \ol{\G}^i_{jk}=0$ (Remark \ref{Rem2}) plus formula
\eqref{R4}. $\square$

 Now, we can easily prove the following statement:
\begin{prop}\label{prop20}
 Let $(M,g)$ be a Brinkmann space with a fixed Brinkmann chart $\{u,v,x^i\}$. If $(M,g)$ is a $2$nd-symmetric manifold, the sections $\wh{B}$ and $\wt{R}$ on $\M$ vanish.
\end{prop}

{\it Proof}.
By $2$nd-symmetry, the integrability conditions associated to the first in (\ref{v21}) and in (\ref{v23}) read, respectively:
\beq\label{ce1} \oln_{[n}\oln_{m]}\wh{B}_{ijkl}=0\Longrightarrow {\ol{R}^r}_{snm}\wh{B}_{rijk}+{\ol{R}^r}_{inm}\wh{B}_{srjk}+{\ol{R}^r}_{jnm}\wh{B}_{sirk}+{\ol{R}^r}_{knm}\wh{B}_{sijr}=0,\\
\label{ce2}\oln_{[n}\oln_{m]}\wh{A}_{sij}=B^r~_{mn}\wh{B}_{s(ij)r}\Longrightarrow \ol{R}^r~_{snm}\wh{A}_{rij}+\ol{R}^r~_{inm}\wh{A}_{srj}+\ol{R}^r~_{jnm}\wh{A}_{sir}=2 B^r~_{nm}\wh{B}_{s(ij)r} .\eeq
In both cases we used \eqref{Riciden}, and for the last case we also used the first in (\ref{v21}) and in \eqref{v23}, as well as \eqref{BBB}. 
If we differentiate \eqref{ce2} with respect to $\oln_k$, using the same information as before plus \eqref{se11}, \eqref{q4} and \eqref{ce1}, we obtain 
$$
\wh{B}_{s(ij)}{}^r \wh{B}_{lrnm}=0 .
$$
Therefore, $\wh{B}$ satisfies all the hypotheses in Corollary \ref{lm1} at each leaf of $\M$, so that $\wh{B}=0$ and, by the identity (\ref{b1}), $\wt{R}=0$ too.
$\square$

%
%

This actually implies, due to \eqref{sss11}, \eqref{se11}, \eqref{v24} and \eqref{v23}, that the three sections $\ol R$, $\wt B$ and $\wh A$ are $\oln$-parallel and $D_0$-parallel.

The integrability equations associated to (\ref{se2}-\ref{se11}), \eqref{qe2}, \eqref{q4} and the first in \eqref{v24} simplified via using Proposition \ref{prop20} are written, on using \eqref{BBB} for the expressions containing $\oln t$, as:
\begin{align}
\nonumber&2\oln_{[n}\oln_{m]}\wt{B}_{ijk}=0,& &2\oln_{[n}\oln_{m]}\wt{A}_{ij}=B^r~_{mn}(2\wt{B}_{(ij)r}-\wh{A}_{rij}),& &\oln_{[n}\oln_{m]}\ol{R}_{ijkl}=0,\\
\nonumber&2 \oln_{[n}\oln_{m]}B_{ijk}=B^r~_{mn}\ol{R}_{r(ij)k},& &\oln_{[n}\oln_{m]}{A}_{ij}= B^r~_{mn}{B}_{(ij)r},&
\end{align}
which via \eqref{Riciden} provide
\beq
\label{ce0}\ol{R}^r~_{inm}\wt{B}_{rjk}+\ol{R}^r~_{jnm}\wt{B}_{irk}+\ol{R}^r~_{knm}\wt{B}_{ijr}&=&0,\\
\label{ce3}\ol{R}^r~_{inm} \wt{A}_{rj}+\ol{R}^r~_{jnm} \wt{A}_{ir}&=&B^r~_{nm}(2\wt{B}_{(ij)r}-\wh{A}_{rij}),\\
\label{cee1}{\ol{R}^r}_{snm}\ol{R}_{rijk}+{\ol{R}^r}_{inm}\ol{R}_{srjk}+{\ol{R}^r}_{jnm}\ol{R}_{sirk}+{\ol{R}^r}_{knm}\ol{R}_{sijr}&=&0,\\
\label{ce5}\ol{R}^r~_{inm} B_{rjk}+\ol{R}^r~_{jnm} B_{irk}+\ol{R}^r~_{knm} B_{ijr}&=&B^r~_{nm}\ol{R}_{rijk},\\
\label{ce4}\ol{R}^r~_{inm}{A}_{rj}+\ol{R}^r~_{jnm}{A}_{ir}&=&2  B^r~_{nm}{B}_{(ij)r}.
\eeq
On the other hand, the integrability conditions derived from \eqref{conmF} 
applied to $\wh{A}$, $B$ and $\wt{A}$ become
\beq\label{ce8}B_{ms}~^r\wh{A}_{rij}+B_{mi}~^r\wh{A}_{irs}+B_{mj}~^r \wh{A}_{ijr}&=&0,\\
\label{ce6}{B_{rjk}B_{mi}~^r+B_{irk}B_{mj}~^r+B_{ijr}B_{mk}~^r}&=&\ol{R}^r~_{ijk}A_{rm},\\
\label{ce7}B_{mi}~^r\wt{A}_{rj}+B_{mj}~^r\wt{A}_{ir}&=&(2\wt{B}_{(ij)r}-\wh{A}_{rij})A^r~_m.
\eeq
where we have used that $\wt B$ and $\wh A$ are $D_0$- and $\oln$-parallel together with \eqref{AAA}, (\ref{se1}-\ref{se2}) for $\wt{A}$ and (\ref{q4}) reduced via Proposition \ref{prop20} for $B$.

Using these integrability equations we can prove
the following further improvement of Proposition \ref{prop20}, which then completes the proof of Theorem \ref{redu}.
\begin{prop}\label{prop21}
Under the hypotheses of Proposition \ref{prop20} the sections $\wt{B}$ and $\wh{A}$ vanish.
\end{prop}
{\it Proof:} We use the results of Proposition \ref{prop20} throughout the proof to simplify the formulas to be combined in the calculations. Start by
applying $D_0$ to \eqref{ce5}, and use \eqref{sss11}, \eqref{qe3}, \eqref{cee1} and \eqref{ce0} in order to get
\beq\label{ee1}\wt{B}^r~_{mk}\ol{R}_{sijr}=0.\eeq
Apply $\oln_n$ to \eqref{ce6} and use \eqref{se11}, \eqref{q4}, \eqref{qe2} and \eqref{ce5} to obtain
\beq\label{ee2}\wh{A}^r~_{mk}\ol{R}_{sijr}=0.\eeq
Applying $D_0$ to \eqref{ce3} and using (\ref{ce2}), the second in \eqref{v23} and \eqref{v24}, \eqref{se1}, \eqref{sss11}, \eqref{qe3} and \eqref{ce0} we get
\beq\label{-}\wt{B}^r~_{mk}(2\wt{B}_{(ij)r}-\wh{A}_{rij})=0.\eeq
However, if we apply $D_0$ to \eqref{ce4} and use  \eqref{sss11}, \eqref{qe3},  \eqref{qe1}, \eqref{ce3} and \eqref{ce5} we derive
$ 2\wt{B}^r~_{mk}B_{(ij)r}+B^r~_{mk}\wh{A}_{rij}=0.$ By applying $D_0$ once more to this last expression and using the second in \eqref{v23} and \eqref{v24}, \eqref{qe3}, \eqref{ee1} and \eqref{ee2} we also get
\beq\label{+}\wt{B}^r~_{mk}(2\wt{B}_{(ij)r}+\wh{A}_{rij})=0.\eeq
In conclusion, comparing \eqref{-} and \eqref{+}, we have that $$\wt{B}^r~_{mk}\wt{B}_{(ij)r}=0$$
so that $\wt{B}$ satisfies all the hypotheses in Proposition \ref{lm2} for each leaf of $\M$, and thus $\wt{B}=0$. By the identity \eqref{b2}, we also get $\wh{A}_{[ij]k}=0$, i.e., $\wh{A}_{ijk}=\wh{A}_{(ij)k}$. $D_0$-differentiating \eqref{ce7}, using \eqref{ce8} and putting $\wt{B}=0$ in  \eqref{se1}, \eqref{qe3} and \eqref{ce3}, we get $\wh{A}_{nrk}\wh{A}^r~_{ij}=0.$ Contracting all the indices and using that $\wh{A}$ is symmetric in the first two, we arrive at $0=\wh{A}_{nrk}\wh{A}^{rnk}=\wh{A}_{rnk}\wh{A}^{rnk}$, i.e., $\wh{A}=0$.
 $\square$

 Therefore, 
\begin{coro}\label{otro}
$\wt A$ is also a $D_0$- and $\oln$-parallel section.
\end{coro}
\begin{rema}\label{tsns}
{\em 
We simply remark that all the equations written in this section, be them 2nd-symmetry ones or
their integrability conditions, are satisfied for any Brinkmann chart: if we perform a new decomposition of type \eqref{chcoo}, the aforementioned equations must be satisfied in the new Brinkmann chart too.

}
\end{rema}

\subsection{Transformation into two independent Lorentzian problems}\label{s8c}

\subsubsection{Reducibility of $\ol{g}$, $\ol{Ric}$ and $\wt{A}$}
For the following application of the generalized Eisenhart Theorem recall Definitions \ref{foli} and \ref{decg}.
  
\begin{prop}\label{ri}
 Let $(M,g)$ be a $2$nd-symmetric Brinkmann space. Then, any Brinkmann decomposition $\{u,v\}$ 
is spatially reducible $$\M=\M^{(1)}\times\M^{(2)}, \hspace{1cm} \ol{g}=\ol{g}^{(1)}\oplus\ol{g}^{(2)}$$ where $\M^{(1)}$ is a $d$-dimensional flat foliation and $\M^{(2)}$  is a $d'$-dimensional locally symmetric foliation with $d,d'\geq 0$, $d+d'=n-2$. Furthermore, $\M^{(2)}$ is itself the product of $s_2$ locally symmetric (non-Ricci flat) Einstein foliations: $$\M^{(2)}=\M^{(2,1)}\times\ldots\times\M^{(2,s_2)}\hspace{1cm}  \ol{g}^{(2)}=\ol{g}^{(2,1)}\oplus\ldots\oplus\ol{g}^{(2,s_2)}$$. \end{prop}

{\it Proof}. As $\ol R$ is $\oln$- and $D_0$-parallel it follows that $D_0 \ol{Ric}=0$ and $\oln~ \ol{Ric}=0$. Then, if the foliation is $u$-Einstein, the result is trivial (put $\M=\M^{(1)}$, $\ol{g}=\ol{g}^{(1)}$ for the Ricci-flat case; otherwise, from $D_0 \ol{Ric}=0$ the foliation is Einstein and $\M=\M^{(2)}$, $\ol{g}=\ol{g}^{(2)}$). If it is not $u$-Einstein, apply Theorem \ref{gei} to $\ol{L}=\ol{Ric}$ and, if $\ol{Ric}$ happens to have a vanishing eigenvalue, choose ${\M}^{(1)}$ as the corresponding Ricci-flat part which, by virtue of Proposition \ref{redfol}(2), it is actually flat. $\square$

Therefore, under the conditions of Proposition \ref{ri} with $d,d'>0$, the partition of the indices corresponds to the reducibility of $\ol{g}$ according to Definition \ref{decg}, where $I_1$ yields the indices of the flat foliation and $I_2$ the indices of the locally symmetric one.  Moreover, the associated curvature tensors satisfy $\ol{R}^{(1)}=0$ and  $\ol{R}^{(2)}\neq 0$, so that $\ol{R}=\ol{R}^{(2)}$ and $\oln ~\ol{R}^{(2)}=0$. Besides, $$\ol{Ric}=\sum_{m=1}^{s_2} \mu_{m}\ol{g}^{(2,m)}, \hspace{1cm} \mu_m\in\R-\{0\}$$ and the flat coordinates $\{x^2,\ldots,x^{d+1}\}$ correspond to the zero eigenvalue $\mu_0=0$. 
\begin{conve}\label{conve}
 {\rm From now on:

(1) When dealing with a \tsns Brinkmann manifold and using Proposition \ref{ri}, we will restrict ourselves to spatially reducible Brinkmann decompositions $\{u,v\}$ so that $\M=\M^{(1)}\times\M^{(2)}$ as in the proposition. This includes the limit case when $\M$ is Einstein, so that either $\M=\M^{(1)}$ ($d'=0$) or $\M=\M^{(2)}$ ($d=0$).

(2) The indices $a,b,c,\ldots$ will run from $2$ to $d+1$ and the indices $a',b',c',\ldots$ will run from $d+2$ to $n-1$. 

(3) We will denote by $\mu^*$ any of the {\it non-zero} eigenvalues of $\ol{Ric}$.
}
\end{conve}

Next, our aim is to prove that $\wt{A}$ is reducible and it admits a similar decomposition to $\ol{Ric}$. Nevertheless, as a difference with $\ol{Ric}$, the non-trivial part of $\wt{A}$ lies on $\M^{(1)}$. 

\begin{prop}\label{teo1}
Choose any spatially-reducible Brinkmann decomposition of a \tsns Brinkmann space. Then, $\wt{A}$ is reducible as $\wt{A}^{(1)}\oplus\wt{A}^{(2)}$ with $\wt{A}^{(2)}=0$ . In addition, there exists a Brinkmann chart $\{u,v,x^i\}$ such that $\ol{g}^{(1)}=\delta_{ab} dx^{a} dx^b$ and the matrix of components $(\wt{A}^{(1)}_{ij})$ is constant and diagonal. 
\end{prop}

{\it Proof}. 
Applying $D_0$ to \eqref{ce6} and using \eqref{qe1}, \eqref{qe3} and Proposition \ref{prop21} together with the fact that $\ol R$ is $D_0$-parallel one derives
$$\ol{R}^r{}_{ijk} \wt{A}_{rm}=0$$
from where we get $\ol{R}^r~_i \wt{A}_{rj}=0$ for all $i,j.$ However, $\wt{A}\neq 0$ by Corollary \ref{Aten} so that, in the Brinkmann chart associated to the spatial reduction of Proposition \ref{ri}, if we set $i\in I_2$, i.e., $i=b'$, we get that $\mu^*\wt{A}_{b'j}=0$, that is, $\wt{A}_{b'j}=0$ for all $j$ and $$\wt{A}=\wt{A}(u,x^i)_{bc} \ol{d}x^{b} \ol{d}x^{c}.$$
Now, from Corollary \ref{otro} we have $0=\oln_{a'} \wt{A}_{bc}=\p_{a'} \wt{A}_{bc}$ (last equality because $\g^a_{b a'}=0$), that is, $\wt{A}$ is reducible with $\wt{A}^{(2)}=0$ and satisfies the hypotheses of Theorem \ref{gei}. Moreover, the coordinate vector fields $\{\p_{d+2},\ldots,\p_{n-1}\}$ span a subspace $S$ of the eigenspace $S_0$ associated to the zero eigenvalue  of $\wt{A}$. Therefore, one can follow the steps of the proof of Theorem \ref{gei}, but working just on $S^\perp$, 
($X(x^{a})=0,$ for all $X\in S$) to obtain
\beit
\item $\ol{g}^{(1)}=\ol{g}^{(1,1)}\oplus\ldots\oplus\ol{g}^{(1,s_1)}.$
\item $\wt{A}=\sum_{m=1}^{s_1} \la_m \ol{g}^{(1,m)}.$
\eit 
\ni As $\ol{g}^{(1)}$ is a flat metric, applying Proposition \ref{redfol}(3) to each of the mutually orthogonal $\ol{g}^{(1,m)}$  a further change of coordinates in the flat block yields $\ol{g}^{(1)}=\delta_{ab} dx'^{a} dx'^b$ (for example, use the proof of that proposition on the restriction to the block of coordinates associated to each $\M^{(1,m)}$). In conclusion, in this Brinkmann chart the matrix of the tensor field $\wt{A}$ 
is a diagonal matrix of constants. $\square$

\begin{rema}\label{410}{\rm Observe that the Brinkmann chart  obtained in Proposition \ref{teo1}:

(1) 
 maintains  the Brinkmann decomposition $\{u,v\}$ of Proposition \ref{ri}. As the tensor fields $\ol{Ric}$ and $\wt{A}$ (Corollary \ref{Aten}) depend only on the decomposition, the conclusions obtained for
 them are independent on the remainder coordinates of the Brinkmann chart.

(2) is such that $\wt{A}$ and $\ol{Ric}$ are orthogonal and occasionally they may vanish simultaneously for a subundle on $\M^{(1)}$.

(3) has a flat metric  $\ol{g}^{(1)}=\delta_{ab} dx^{a} dx^b$ so that $\oln_{a}=\p_{a}$. Moreover, $\oln_{a'} T^{a}_{b}=\p_{a'} (T^{a}_{b}).$
}
\end{rema}

\subsubsection{The building blocks of the metric $g$}\label{s8b}

Until now, we have proven the reducibility of the metric $\ol{g}$ for any Brinkmann decomposition $\{u,v\}$. Our aim is to prove that, in fact, the $2$nd-symmetry induces the existence of two simpler $2$nd-symmetric Brinkmann spaces associated to the original Brinkmann space $(M,g)$. 
\begin{lema}\label{lem3}
For a Brinkmann chart $\{u,v,x^i\}$ as in Proposition \ref{teo1}, the sections $B$ and $A$ in Definition \ref{AB} satisfy that:
\begin{center}
\begin{tabular}{lllll}
 $B_{{a}{b}a'}=0$;&$B_{a'{a}{b}}=0$;&$B_{a'b'{a}}=0$;& $B_{{a}a'b'}=0$;&
$A_{{a}a'}=0.$
\end{tabular}
\end{center}
\end{lema}

{\it Proof}.  
The commutation property \eqref{conmF} applied to the $\oln$- and $D_0$ parallel section $\ol R$ provides the following integrability condition
$$
B_{ms}~^r\ol{R}_{rijk}+B_{mi}~^r\ol{R}_{srjk}+B_{mj}~^r\ol{R}_{sirk}+B_{mk}~^r\ol{R}_{sijr}=0.
$$
(i) Contracting here $i$ and $k$ and
\begin{itemize}
 \item taking $m={a}$, $s={b}$, $j=a'$, one gets:  $B_{{a}{b}}~^r \ol{R}_{r{a'}}=B_{{a}{b}a'}~\mu^*=0$ that is, $ B_{{a}{b}a'}=0,$
\item taking $m=a'$, $s=b'$, $j={a}$, one gets: $B_{a'{a}}~^r \ol{R}_{b'r}=B_{a'{a}b'}~\mu^*=0$ that is, $B_{a'b'{a}}=0.$
\end{itemize}
(ii) Using now $B_{[ijk]}=0$, and $B_{i[jk]}=0$, (i) above and:
\begin{itemize}
 \item taking $i=a'$, $j=b'$ and  $k=a$, one has that $B_{aa'b'}=0$,
\item taking $i=a'$, $j=a$ and  $k=b$, one has that $B_{a'ab}=0$.
\end{itemize}
(iii) Contracting $i$ and $s$ in the equation (\ref{ce4}) and
taking $j=a, m=a'$, one derives $\ol{R}_{ra'}A^r~_{a}=2B_{(ia)}~^rB_{ra'}~^i$. Using here (i) and (ii) above one arrives at 
$\mu^* A_{a'a}=0\Longrightarrow A_{aa'}=0. ~\square$

\begin{prop}\label{cor1}
 For a Brinkmann chart  $\{u,v,x^i\}$ as in Proposition \ref{teo1}, the sections $t$ and $h$ 
 are reducible.
 \end{prop}

{\it Proof}. We have to prove that
\begin{center}
\begin{tabular}{llllll}
 $t_{{a}a'}=0$;&$t_{a'{a}}=0$;&$\p_{a'} t_{{a}{b}}=0$;& $\p_{a} t_{a'b'}=0$;&$\p_{a} h_{a'}=0$;&$\p_{a'} h_{a}=0$
\end{tabular}
\end{center}
The fact that $\ol R$ is $D_0$-parallel can be rewritten using \eqref{D0T} as
$$\dot{\ol{R}}_{ijkl}+t^r~_i \ol{R}_{rjkl}+t^r~_j \ol{R}_{irkl}+t^r~_k \ol{R}_{ijrl}+t^r~_l \ol{R}_{ijkr}=0.$$ 
Taking $i=a$ one gets that $t^r~_{a} \ol{R}_{rjkl}=0$, so that contracting $j$ and $l$ and taking $k=a'$ one gets 
\beq\label{tab}0=t_{a'a}=-t_{aa'}.\eeq 
By \eqref{BBB} and Lemma \ref{lem3} it follows that $0=B_{{a}{b}a'}=\oln_{a'} t_{{a}{b}}-\oln_{b} t_{{a}a'}$, which implies by Remark \ref{410}(2) and
\eqref{tab} that $\partial_{a'} t_{{a}{b}}=0$. Analogously, $0=B_{a'b'{a}}=\oln_{a} t_{a'b'}-\oln_{b'} t_{a'a}$ implies $\partial_{{a}} t_{a'b'}=0$. Lemma \ref{lem3}  and \eqref{AAA} yield $0=A_{{a}a'}=t_{{a}r}t^r~_{a'}-\oln_{a'} h_{a}-D_0 t_{{a}a'}$, so sustituting \eqref{tab} and using \eqref{D0T} and Remark \ref{410}(2) again one obtains $\partial_{{a}} h_{a'}=\partial_{a'} h_{{a}}=0$.
$\square$

\begin{teor}\label{cor2}
 For a \tsns Brinkmann space, there exists a spatially reducible Brinkmann decomposition $\{u',v'\}$ such the function $H$ and the one-form section $W$ in (\ref{m1}) are reducible in the associated Brinkmann chart $\{u',v',x'^i\}$. This chart is related to that obtained in Proposition \ref{teo1} by $\{u',v',x'^i\}=\{u,v+f(u,x^j),x^i\}$ for some function $f$.
\end{teor}

{\it Proof}. We must prove that there exists a Brinkmann chart  $\{u',v',x'^i\}$ such that
$$H(u',x'^i)=H^{(1)}(u',x'^{a})+ H^{(2)}(u',x'^{a'}),$$
$$W_{a}=W_{a}(u',x'^b)(=W^{(1)}_{a}(u',x'^b)),~~~~
W_{a'}=W_{a'}(u',x'^{b'})(=W^{(2)}_{a'}(u',x'^{b'})).$$
Let  $\{u,v,x^i\}$ be the Brinkmann chart obtained in Proposition \ref{teo1}, so that  by Proposition \ref{cor1} $t$ and $h$ are reducible. 

Simplification of $W$: $\partial_{a'} t_{{a}{b}}=0$ implies that $W_{{a},{b}}-W_{{b},{a}}$ depends only on the coordinates $\{u,x^{{a}}\}$, hence there exists a function $f_1(u,x^i)$ such that $W_{{a}}(u,x^{j})=f_{1,{a}}(u,x^j)+w_{{a}}(u,x^{{b}}).$ Analogously, $\partial_{a} t_{{a'}{b'}}=0$ implies the existence of
a function $f_2(u,x^i)$ such that $W_{a'}(u,x^i)=f_{2,a'}(u,x^i)+w_{a'}(u,x^{b'})$.
Then, \eqref{tab} tells us that $f_{1,a,a'}=f_{2,a',a}$ so that $f_{1}(u,x^i)-f_{2}(u,x^i)=F(u,x^a)+G(u,x^{a'})$ for some functions $F$ and $G$ which can be absorbed into $w_a(u,x^b)$ and $w_{a'}(u,x^{b'})$, respectively. Consequently, there exists a function $f(u,x^j)$ such that
\beq\label{ww} W_{a}(u,x^b)=f_{,{a}}(u,x^j)+w_{a}(u,x^{b}),~~~~
W_{a'}(u,x^{b'})=f_{,a'}(u,x^j)+w_{a'}(u,x^{b'}).
\eeq
Simplification of $H$: using \eqref{h}, $\partial_{{a}} h_{a'}=\partial_{a'} h_{{a}}=0$ provide
$H_{,{a}} (u,x^i)=h_{{a}} (u,x^{{b}})+\dot{W}_{{a}}(u,x^{i})$ and 
$ H_{,a'} (u,x^i)=h_{a'} (u,x^{b'})+\dot{W}_{a'}(u,x^i)$, which become after use of  \eqref{ww},
$H_{,{a}} (u,x^i)=h_{{a}} (u,x^{{b}})+\dot{f}_{,{a}}(u,x^i)+\dot{w}_{{a}}(u,x^{{b}})$ and $H_{,a'} (u,x^i)=h_{a'} (u,x^{b'})+\dot{f}_{,a'}(u,x^i)+\dot{w}_{a'}(u,x^{b'}),$ from where it is easy to deduce that
\beq\nonumber
H(u,x^i)=\dot{f}(u,x^i)+H^{(1)}(u,x^{{a}})+ H^{(2)}(u,x^{a'})\eeq
\noindent for some functions $H^{(1)},H^{(2)}$.

By choosing now the new Brinkmann decomposition defined by $v'=v+f(u,x^i)$, the thesis follows on using \eqref{H} and \eqref{Wi}.
$\square$


\begin{rema}\label{rm1}{\rm
(1)  What we have proven is that there exists a Brinkmann chart $\{u,v,x^i\}$ and a partition of the indices $I_1=\{2,\dots, d+1\}, I_2=\{d+2,\dots, n-1\}$ for some $d\in
\{0,\dots , n-2\}$ (recall Convention \ref{conve}), such that $\ol{g}$, $H$ and $W$ are simultaneously reducible in the sense given of Definition \ref{decg}.

(2) Recall that if $\{u,v\}$ is spatially reducible, there exist two foliations $\M^{(1)},\M^{(2)}$ with associated leaves $(\ol{M}^{(1)},\ol{g}^{(1)}),(\ol{M}^{(2)},\ol{g}^{(2)})$ such that $\M=\M^{(1)}\times\M^{(2)}$, $\ol{M}=\ol{M}^{(1)}\times\ol{M}^{(2)}$ and $\ol{g}=\ol{g}^{(1)}\oplus \ol{g}^{(2)}$ (Remark \ref{*}). Observe then that the metric $g$ can be written as
$$g=-2 du(dv+(H^{(1)}+H^{(2)})du+\ri W^{(1)}+\ri W^{(2)})+\ri {\ol{g}}^{(1)}\oplus\ri {\ol{g}}^{(2)},$$
so that to any 2nd-symmetric Brinkmann space we can associate a pair 
of lower-dimensional Brinkmann spaces $(M^{[m]},g^{[m]}), m\in\{1,2\}$ by
$M^{[m]}=\R^2\times \ol{M}^{(m)}$  and
$$
g^{[m]}= -2du(dv + H^{(m)}du+ W^{(m)}) + \ri {\ol{g}}^{(m)}.
$$
Of course, such a pair of spaces may be non-unique. $(M^{[m]},g^{[m]})$ are the building blocks of any proper 2nd-symmetric Lorentzin manifold, and they are extremely helpful in the resolution of the equations for 2nd-symmetry because they actually simplify to the equations corresponding to each of
the two simpler Brinkmann spaces $(M^{[m]},g^{[m]}), m\in\{1,2\}$ as the next proposition proves.}
\end{rema}

\begin{prop}\label{redBri}
The pair $(M^{[m]},g^{[m]})$ of Brinkmann
spaces associated to any 2nd-symmetric Brinkmann space  $(M,g)$ according the the previous Remark are, themselves, 2nd-symmetric.
\end{prop}
{\em Proof.}  Note first of all that $\ol{R}, h$ and $t$ are reducible, so that $A$, $B$ and their $\oln$- and $D_0$-derivations are reducible too, from  where $\wt{A}$ is also reducible. Therefore, for $m\in\{1,2\}$, the reduced sections $\ol{R}^{(m)}$, $h^{(m)}$, $t^{(m)}$, $A^{(m)}$, $B^{(m)}$ correspond to the geometrical objects for the Brinkmann spaces $(M^{[m]},g^{[m]})$ as defined in \eqref{olRR}, \eqref{h}, \eqref{t} and Definition \ref{AB} respectively.
It is then straightforward to check that, in the Brinkmann chart of Theorem \ref{cor2},
the equations of $2$nd-symmetry for $(M,g)$ are exactly identical with the combination of the equations of $2$nd-symmetry for the two building blocks $(M^{[m]},g^{[m]})$. $\square$

 \begin{rema}\label{r519}{\rm
%

In conclusion, applying Proposition \ref{redBri} to the Brinkmann decompositions $\{u',v'\}$ of Theorem \ref{cor2},  we can reorganize the equations of $2$nd-symmetry in two simpler sets: 
\begin{itemize}
 \item The equations associated to the Brinkmann space $(M^{[1]},g^{[1]})$ with coordinates $\{u,v,x^{a}\}$, such that: $\ol{g}^{(1)}=\delta_{ab} dx^{a} dx^b$, $\ol{R}^{(1)}=0$ (ergo $\ol{\text{Ric}}^{(1)}=0$) and $\wt{A}^{(1)}=\sum_{l=1}^{s_1} \la_l \ol{g}^{(1,l)}$ with some $\la_l\neq 0$.
\item Those associated to $(M^{[2]},g^{[2]})$ with coordinates $\{u,v,x^{a'}\}$, such that: $\ol{g}^{(2)}=\ol{g}^{(2,1)}\oplus\ldots\oplus\ol{g}^{(2,s_2)}$, $\oln ~\ol{R}^{(2,l)}=0$ but $\ol{R}^{(2,l)}\neq 0$, $\ol{\text{Ric}}^{(2)}=\sum_{l=1}^{s_2} \mu_{l}\ol{g}^{(2,l)}$ with each $\mu_{l}\neq 0$ and $\wt{A}^{(2)}=0$.
\end{itemize}}
\end{rema}

\subsection{Proof of the main Theorem \ref{lem4}}\label{s11}
As a consequence of Corollary \ref{lnuevo}, we can apply Lemma
\ref{glob} and, consequently, only the local version of the result must be
proven. We will start the computations in the Brinkmann chart
$\{u',v',x'^i\}$ of Theorem \ref{cor2} but dropping the primes for
the sake of clarity in the notation.

By Remark \ref{r519}, we consider first the equations for
$(M^{[2]},g^{[2]})$. Since $\wt{A}^{(2)}=0$,
Corollary \ref{Aten} informs us that this is a locally symmetric Lorentzian manifold with a parallel lightlike vector field.
Therefore, Theorem \ref{CW} implies that it is locally isometric to the
product of a symmetric (non-flat) Riemannian space with
$\text{$\mathbb{L}^2$}=(\R^2,-2dudv)$.  In particular, up to a change of coordinates of
type \beq\label{ch1} u'=u,\hspace{0.5cm}
v'=v+F(u,x^{a'});\hspace{0.5cm} y^{a'}=y^{a'}(u,x^{b'})\eeq
$(M^{[2]},g^{[2]})$ has $H^{(2)}=0$, $W^{(2)}=0$ and $\dot{\ol{g}}^{(2)}=0$.

Consider now the $2$nd-symmetry equations associated to
$(M^{[1]},g^{[1]})$ according to Remark \ref{r519}. Using $\ol
R^{(1)}=0$  in  (\ref{ce4})
we have that $B_{({a}{b})c}B^{c}~_{de}=0$, hence Proposition \ref{lm2} leads us to
\beq\label{b0} B^{(1)}=0.\eeq
From \eqref{BBB} follows that
$\oln_{c} t_{{a}{b}}-\oln_{b} t_{{a}{c}}=0,$ which together with \eqref{t} and Remark \ref{410}(2) provides $\p_{c} t_{{a}{b}}-\p_{b} t_{{a}{c}}=\p_{a} t_{{c}{b}}=0.$ We conclude that $t^{(1)}$ is $\oln$-parallel. We have two consequences of this fact:
\beit
\item $W_{{a},{b}}-W_{{b},{a}}=2t_{ab}(u)$ depend only on $u$. Consequently, $W^{(1)}$ can be written as $$W_{a}(u,x^{b})=f_{,{a}}(u,x^{b})+t_{{a}{c}}(u)x^{c}$$  for some function $f(u,x^{b})$.
\item Setting $\wh A=0$ in \eqref{qe2} and by Lemma \ref{lem3} and \eqref{b0} we have that $\p_{c} A_{{a}{b}}=0$, which from \eqref{AAA} and Proposition \ref{cor1} implies that $\p_{c}\p_{b} h_{a}=0$, that is, $h^{(1)}$ takes the form $$h_{a}(u,x^{b})=\Lambda_{{a}{c}}(u)x^{c}+B_{a}(u)$$ for some functions $\Lambda_{ac}$ and $B_a$. From \eqref{h} we then have
\beq \nonumber H_{{,a}}(u,x^{b})&=&\Lambda_{{a}{c}}(u)x^{c}+B_{a}(u)+\dot{W}_{a}(u,x^{b})\\
\nonumber &=&\Lambda_{{a}{c}}(u)x^{c}+B_{a}(u)+\dot{f}_{,{a}}(u,x^{b})+\dot{t}_{{a}{c}}(u)x^{c}.
\eeq
Observe that the symmetry of $H_{,a,b}$ implies that $\Lambda_{ab}+\dot{t}_{ab}=\Lambda_{ba}+\dot{t}_{ba}=\Lambda_{(ab)}$ where the antisymmetric character of $\dot{t}$ has been used in the last equality.
\eit
In conclusion, the metric for $M^{[1]}$ becomes
\beq \label{ch4} g^{[1]}&=&-2du\left(dv+H^{(1)}(u,x^{a}) du+\ri W^{(1)}\right)+\delta_{ab} dx^{a}dx^{b}\eeq
\ni with
\beq \nonumber W^{(1)}&=&W_{a}(u,x^{b})\ol{d}x^a =\left[f_{,{a}}(u,x^{b})+t_{{a}{c}}(u)x^{c}\right]\ol{d}x^a;\\
\nonumber
H^{(1)}(u,x^{a})&=&\dot{f}(u,x^{a})+\frac{1}{2}\Lambda_{({b}{c})}(u)x^{b} x^{c}+ B_{c}(u)x^{c}+C(u).
\eeq

Next, we use the following claim, to be proven later.
\begin{claim}
 For $(M^{[1]},g^{[1]})$, there exists a change of Brinkmann chart of type
\beq\label{ch2}
 u'=u,~~~v'=v+\chi(u,x^{a}),~~~y^{a}= R^a_b(u) x^{b} + D^{a}(u),\eeq
\noindent such that the metric becomes
\beq\label{ch3}
  g^{[1]}=-2du'\left(dv'+H'(u',y^{a}) du'\right)+\delta_{ab}dy^{a}dy^{b}.
\eeq
\ni where $H'(u',y^{a})=-A_{{a}{b}}(u')y^{a} y^{b}$, $A_{ab}$ being the components of the section $A$ in Definition \ref{AB} associated to the Brinkmann chart $\{u',v',y^i\}$.
\end{claim}

Combining the change of coordinates in this claim and the one given in \eqref{ch1}, there exists a change of coordinates in the  entire Brinkmann space $(M,g)$ of type
\begin{align}
 \nonumber & u'=u, & &  v'=v+F(u,x^{a'})+\chi(u,x^{a}), \\
 \nonumber & y^{a'}=y^{a'}(u,x^{b'}),& & y^{a}= R^a_b(u) x^{b} + D^{a}(u)
\end{align}
such that the metric of $(M,g)$ becomes
 $$ g=-2du'\left(dv'+H'(u',y^{a}) du'\right)+\delta_{ab} dy^{a}dy^{b}+\ol{g}_{a'b'}(y^{c'})dy^{a'} dy^{b'}$$
where $H'(u',y^{a})=-A_{{a}{b}}(u')y^{a} y^{b}$. 
 To end the proof, note that $\wt A=\dot A$ and therefore Corollay \ref{otro} gives, via $D_0\wt A =0$, $\ddot{A}_{{a}{b}}(u)=0.$
Of course, we need $\dot{A}_{{a}_0{b}_0}(u)\neq 0$ for some ${a}_0,{b}_0$ in order for the manifold not to be locally symmetric, due to Corollary \ref{Aten}.

{\it Proof of the claim}.
In order to find the required $\chi$, $D^a$ and $R^a_b$, put $H'(u',y^a)=-A_{ab}(u')y^a y^b$, substitute \eqref{ch2} in \eqref{ch3} and require that the obtained expression equals \eqref{ch4}. Then, the following equations arise:
\begin{align}
\nonumber & \dot{\chi}=H^{(1)}-H'+\fr 1 2 \delta_{ab}(\dot{R}^a_c x^c+\dot{D}^a)(\dot{R}^b_d x^d+\dot{D}^b),\\
\nonumber &\chi_{,a}=W_a+\fr 1 2\delta_{bc}\left((\dot{R}^b_d x^d+\dot{D}^b)R^c_a+(\dot{R}^c_dx^d+\dot{D}^c)R^b_a\right),\\
\label{orto} &\delta_{ab} R^a_c R^b_d=\delta_{cd}.
\end{align}
Here, the known data are $t_{ab}(u)$, $B_{c}(u)$ and $\Lambda_{(ab)}(u)$, while the unknowns are $R^{a}_{b}(u)$, $D^b(u)$ and $\chi(u,x^{a})$.
The integrability conditions (given by the cross derivatives) of the first two expressions yield
\begin{align}
\label{eqq1} &t_{cd}=\fr 1 2 \delta_{ab}(\dot{R}^a_c R^b_d-R^b_c\dot{R}^a_d),\\
\label{eqq2}&\Lambda_{(cd)}=-2A_{be} R^b_c R^e_d-\fr 1 2 \delta_{ab}(R^a_c \ddot{R}^b_d+R^b_c \ddot{R}^a_d),\\
\label{eqq3}& B_{c}=-2 A_{be}R^b_c D^e+ \delta_{ab}\ddot{D}^a R^b_c.
\end{align}
To prove that these equations have solutions, proceed as follows.
Define $R^c_{b}(u)$ as a solution of the ODE
$$
\dot R^c_{b}=-\delta^{dc}(R^{-1})^{a}_{d}\, t_{ab}
$$
for the given $t_{ab}(u)$. Then, equation \eqref{eqq1} is
automatically satisfied. Concerning \eqref{orto}, note that, for
such a solution, the derivative of
$\delta_{cd}R^c_{a}R^d_{b}$ 
vanishes due to the anti-symmetry of $t_{ab}(u)$. Therefore, by
imposing as initial condition that $R^b_{a}(0)$ is
any orthogonal matrix, the necessary condition \eqref{orto} holds
for all $u$ (that is, $R^a_b(u)$ is a curve of rotation matrices,
which can be chosen by using $d(d-1)/2$ free parameters
codified in $R^b_{a}(0)$).
Once $R^b_{c}(u)$ is determined, equation \eqref{eqq2} fixes
$A_{ab}(u)$ and, using this, the existence of $D^b(u)$ (which will
depend on two new constant vectors, that is, on $2d$ new parameters) and, a posteriori, of $\chi$, is
ensured by \eqref{eqq3} plus, again, standard results in
differential equations.
 $\square$

\begin{rema}{\rm As $(M_1,g_1)$ in Theorem \ref{lem4} is a proper Cahen-Wallach space of order $2$, let $O^i_j$ be the matrix of rotations that diagonalizes the symmetric matrix $A^{(1)}$ in \eqref{ecw}. Then, performing a change of coordinates of type $u'=u-u_0$ and $y^i=O^i_j x^j$ we can write $A^{(1)}$ as a diagonal matrix and cancel one of the elements of the symmetric matrix $A^{(0)}$. In conclusion, the number of essential parameters of a proper $2$nd-symmetric Lorentzian manifold with fixed $K=-\p_v$ is given by $d-1+d(d+1)/2$ (observe that $A^{(1)}$ is a $d\times d$ square matrix). If $K$ is not fixed, by a change of coordinates of type $u'=\varepsilon u-u_0$ with $v'=v/\varepsilon$ for some constant $\varepsilon\neq 0$ and $y^i=O^i_j x^j$, the same simplifications can be achieved and furthermore one of the non-zero eigenvalues of $A^{(1)}$ can be set to $\pm 1$, in which case the number of essential parameters is one less.}
\end{rema}


\end{document}